\documentclass[a4paper,12pt]{amsart}
\usepackage{amssymb}
\usepackage{ifthen}
\usepackage{graphicx}
\nonstopmode \numberwithin{equation}{section}
\usepackage[usenames]{color}
\newtheorem{thm}{Theorem}
\newtheorem{cor}{Corollary}[section]
\newtheorem{lem}{Lemma}[section]
\newtheorem{prop}{Proposition}
\newtheorem{rem}{Remark}[section]
\newtheorem{rems}{Remarks}
\newtheorem{step}{Step}[subsection]
\newtheorem{defin}{Definition}[section]
\newtheorem{examp}{Example}[section]
\newtheorem{prob}{Problem}
\newtheorem{ques}[equation]{Question}
\newtheorem{op}{Problem}

\newtheorem*{mysolution}{Solution}

\newtheorem{conj}[equation]{Conjecture}
\newtheorem{deter}[equation]{Determination}
\newtheorem{case}{Case}[section]
\newtheorem{subcase}{Subcase}
\newtheorem{subsubcase}{Subsubcase}
\newtheorem{claim}{Claim}[section]
\newtheorem{assertion}{Assertion}[section]
\newtheorem{subclaim}{Subclaim}[subsection]

\newcounter {own}
\def\theown {\thesection       .\arabic{own}}

\newenvironment{pf}[1][]{%
 \vskip 3mm
 \noindent
 \ifthenelse{\equal{#1}{}}%
  {{\bf Proof. }}%
  {{\bf #1.} }%
 }%
{\qed\bigskip}

\newcounter{alphabet}
\newcounter{tmp}
\newenvironment{Thm}[1][]{\refstepcounter{alphabet}%
\bigskip%
\noindent%
{\bf Theorem \Alph{alphabet}}%
\ifthenelse{\equal{#1}{}}{}{ (#1)}%
{\bf .} \itshape}{\vskip 8pt}

\makeatletter
\newcommand{\Ref}[1]{\@ifundefined{r@#1}{}{\setcounter{tmp}{\ref{#1}}\Alph{tmp}}}
\makeatother

\newcommand{\IR}{{\mathbb R}}

\newcommand{\IC}{{\mathbb C}}

\newcommand{\id}{{\operatorname{\mathrm{id}}}}

\newcommand{\diam}{{\operatorname{diam}}}

\newcommand{\dist}{{\operatorname{dist}}}

\newcommand{\bsol}{\begin{mysolution}}
\newcommand{\esol}{\end{mysolution}}

\def\be{\begin{equation}}
\def\ee{\end{equation}}

\newcommand{\ben}{\begin{enumerate}}
\newcommand{\een}{\end{enumerate}}

\newcommand{\blem}{\begin{lem}}
\newcommand{\elem}{\end{lem}}
\newcommand{\bthm}{\begin{thm}}
\newcommand{\ethm}{\end{thm}}
\newcommand{\bcor}{\begin{cor}}
\newcommand{\ecor}{\end{cor}}
\newcommand{\beg}{\begin{examp}}
\newcommand{\eeg}{\end{examp}}
\newcommand{\begs}{\begin{examples}}
\newcommand{\eegs}{\end{examples}}
\newcommand{\bdefe}{\begin{defin}}
\newcommand{\edefe}{\end{defin}}
\newcommand{\bprob}{\begin{prob}}
\newcommand{\eprob}{\end{prob}}
\newcommand{\bques}{\begin{ques}}
\newcommand{\eques}{\end{ques}}
\newcommand{\bei}{\begin{itemize}}
\newcommand{\eei}{\end{itemize}}

\newcommand{\bas}{\begin{assertion}}
\newcommand{\eas}{\end{assertion}}

\newcommand{\bde}{\begin{deter}}
\newcommand{\ede}{\end{deter}}
\newcommand{\bca}{\begin{case}}
\newcommand{\eca}{\end{case}}
\newcommand{\bsca}{\begin{subcase}}
\newcommand{\esca}{\end{subcase}}

\newcommand{\bssca}{\begin{subsubcase}}
\newcommand{\essca}{\end{subsubcase}}

\newcommand{\bcl}{\begin{claim}}
\newcommand{\ecl}{\end{claim}}
\newcommand{\bscl}{\begin{subclaim}}
\newcommand{\escl}{\end{subclaim}}
\newcommand{\bcon}{\begin{conj}}
\newcommand{\econ}{\end{conj}}
\newcommand{\bcons}{\begin{conjs}}
\newcommand{\econs}{\end{conjs}}
\newcommand{\bprop}{\begin{prop}}
\newcommand{\eprop}{\end{prop}}
\newcommand{\br}{\begin{rem}}
\newcommand{\er}{\end{rem}}
\newcommand{\brs}{\begin{rems}}
\newcommand{\ers}{\end{rems}}
\newcommand{\bo}{\begin{obser}}
\newcommand{\eo}{\end{obser}}
\newcommand{\bos}{\begin{obsers}}
\newcommand{\eos}{\end{obsers}}
\newcommand{\bpf}{\begin{pf}}
\newcommand{\epf}{\end{pf}}
\newcommand{\ba}{\begin{array}}
\newcommand{\ea}{\end{array}}
\newcommand{\beq}{\begin{eqnarray}}
\newcommand{\beqq}{\begin{eqnarray*}}
\newcommand{\eeq}{\end{eqnarray}}
\newcommand{\eeqq}{\end{eqnarray*}}

\newcommand{\bop}{\begin{op}}
\newcommand{\eop}{\end{op}}
\newcommand{\bsp}{\begin{step}}
\newcommand{\esp}{\end{step}}

\newtheorem{pfofThm1.5}[equation]{}

\newcounter{minutes}
\divide\time by 60
\newcounter{hours}
\multiply\time by 60 \addtocounter{minutes}{-\time}

\begin{document}

\bibliographystyle{amsplain}

\title{Semisolidity and locally weak quasisymmetry of homeomorphisms in metric spaces}

\def\thefootnote{}
\footnotetext{ \texttt{\tiny File:~\jobname .tex,
          printed: \number\year-\number\month-\number\day,
          \thehours.\ifnum\theminutes<10{0}\fi\theminutes}
} \makeatletter\def\thefootnote{\@arabic\c@footnote}\makeatother

\author{Manzi Huang}
\address{Manzi Huang, Department of Mathematics,
Hunan Normal University, Changsha,  Hunan 410081, People's Republic
of China} \email{mzhuang79@163.com}

\author{Antti Rasila}
\address{Antti Rasila,
 Department of Mathematics, Hunan Normal University, Changsha,  Hunan 410081, People's Republic
of China, and\\
 Department of Mathematics and Systems Analysis, Aalto University, FI-00076 Aalto, Finland}
\email{antti.rasila@iki.fi}


\author{Xiantao Wang$^{\mathbf{*}}$
}
\address{Xiantao Wang, Department of Mathematics, Shantou University, Shantou 515063, People's Republic
of China} \email{xtwang@hunnu.edu.cn}

\author{Qingshan Zhou}
\address{Qingshan Zhou,  Department of Mathematics, Shantou University, Shantou 515063, People's Republic
of China} \email{q476308142@qq.com}

\date{}
 \keywords{
Quasihyperbolic metric; length metric; quasiconformality; quasisymmetry; weak quasisymmetry; locally weak quasisymmetry; semisolidity; quasiconvex metric space.\\
${}^{\mathbf{*}}$ Corresponding author}

\begin{abstract}
In this paper, we investigate the relationship between semisolidity and locally weak quasisymmetry of homeomorphisms in quasiconvex and complete metric spaces. Our main objectives are to $(1)$ generalize the main result in \cite{HL} together with other related results, and $(2)$ give a complete answer to the open problem given in \cite{HL}.
As an application, we prove that the composition of two locally weakly quasisymmetric mappings is a locally weakly quasisymmetric mapping and that it is quasiconformal.
\end{abstract}

\thanks{The research was partly supported by NSF of
China (No. 11571216)}

\maketitle \pagestyle{myheadings} \markboth{Manzi Huang, Antti Rasila, Xiantao Wang and Qingshan Zhou
}
{Semisolidity and locally weak quasisymmetry of homeomorphisms in metric spaces}


\section{Introduction and main results}\label{sec-1}


The quasihyperbolic metric was introduced by Gehring and
his students Palka and Osgood in the 1970's \cite{GO, GP} in the setting of Euclidean spaces
${\mathbb R}^n$ $(n\ge 2).$ Since its first appearance, the quasihyperbolic metric has become
an important tool in the geometric function theory of Euclidean spaces.

From late 1980's onwards, V\"ais\"al\"a has developed the theory of (dimension) free quasiconformal mappings (the free theory) in Banach spaces \cite{Vai-1, Vai-2, Vai-3, Vai-4, Vai-5}, which is based on the quasihyperbolic metric. The main advantage of this approach over generalizations based on the conformal modulus (see \cite{Heinonen-book} and references therein) is that it does not make use of volume integrals, which allows one to study quasiconformality of homeomorphisms in infinite dimensional Banach spaces. In the free theory,
V\"ais\"al\"a has mainly studied the relationships between quasiconformal mappings and quasisymmetric mappings, as well as properties of the quasihyperbolic metric and various classes of domains. The importance of the quasihyperbolic metric in this setting arises from the distortion inequality in Definition \ref{def-1.3} (cf. the Schwarz-Pick type result given in \cite{GO}). This line of research has recently attracted substantial interest in the research community (see e.g. \cite{BB, BHK, BHX, HDA, HSX, hlvw, rt1, rt2, S-W, vai2014}).

The class of quasisymmetric mappings on the real axis were first introduced by Beurling and
Ahlfors \cite{BA}, who found a way to obtain a quasiconformal extension of a quasisymmetric self-mapping of the
real axis to a self-mapping of the upper half-plane. This concept
was later generalized by Tukia and V\"ais\"al\"a, who studied quasisymmetric
mappings between metric spaces \cite{TV}. 

In 1998, Heinonen and Koskela \cite{HK} proved a remarkable result, showing that these two concepts, quasiconformality and
quasisymmetry, are quantitatively equivalent in a large class of metric spaces, which
includes Euclidean spaces. In Banach spaces, V\"ais\"al\"a also proved the quantitative equivalence among free quasiconformality,
quasisymmetry and weak quasisymmetry. See \cite[Theorem 7.15]{Vai-5}.

Another important result in this topic was very recently obtained by Huang and Liu \cite{HL}. They used ideas from V\"ais\"al\"a's approach to show that even weakly quasisymmetric mappings preserve, in a suitable way, the quasihyperbolic metric in quasiconvex and complete metric spaces (see Theorem $1.6$ in \cite{HL}). Also, in \cite{HL}, the authors asked that whether the converse of Theorem $1.6$ holds or not. As the main purpose of this paper, we shall study the main result and the open problem in \cite{HL}. In order to introduce the main result and the problem in \cite{HL} and state our results, we
need some basic definitions.

Throughout this paper, we always assume that $X$ and $Y$ are metric spaces.
Following notations and terminology of \cite{HK-1, HK, HL, Tys, Vai-5}, we begin with the definitions of quasiconformality and quasisymmetry of homeomorphisms.

\bdefe\label{japan-30}
A homeomorphism $f$ from $X$ to $Y$ is said to be \begin{enumerate}
\item
{\it quasiconformal} if there is a constant $H<\infty$ such that
\be\label{mon-1}\limsup_{r\to 0}\frac{L_f(x,r)}{l_f(x,r)}\leq H\ee
for all $x\in X$;
\item
{\it quasisymmetric} if there is a constant $H<\infty$ such that
\be\label{mon-2}\frac{L_f(x,r)}{l_f(x,r)}\leq H\ee
for all $x\in X$ and all $r>0$,\end{enumerate}
 where $L_f(x,r)=\sup_{|y-x|\leq r}\{|f(y)-f(x)|\}$
and
$l_f(x,r)=\inf_{|y-x|\geq r}\{|f(y)-f(x)|\}$.
\edefe

Here and in what follows, we always use $|x-y|$ to denote the distance between $x$ and $y$.

\bdefe\label{japan-31} A homeomorphism $f$ from $X$ to $Y$ is said to be
\begin{enumerate}
\item $\eta$-{\it quasisymmetric} if there is a homeomorphism $\eta : [0,\infty) \to [0,\infty)$ such that
$$ |x-a|\leq t|x-b|\;\; \mbox{implies}\;\;   |f(x)-f(a)| \leq \eta(t)|f(x)-f(b)|$$
for each $t>0$ and for each triplet $x,$ $a$, $b$ of points in $X$;

\item {\it weakly $H$-quasisymmetric} if
$$ |x-a|\leq |x-b|\;\;  \mbox{ implies}\;\;   |f(x)-f(a)| \leq H|f(x)-f(b)|$$
for each triplet $x$, $a$, $b$ of points in $X$.
\end{enumerate}\edefe

\br\label{rem-1.1} The following observations are immediate consequences of Definitions \ref{japan-30} and \ref{japan-31}.
\begin{enumerate}
\item\label{hwz-1}
The quasisymmetry implies the quasiconformality;
\item\label{hwz-2}
A homeomorphism $f$ from $X$ to $Y$ is quasisymmetric with coefficient $H$ defined by \eqref{mon-2} if and only if it is weakly $H$-quasisymmetric;
\item\label{hwz-3}
The quasisymmetry implies the weak quasisymmetry. In general, the converse is not true (cf. \cite[Theorem $8.5$]{Vai-5}). See \cite{hprw}
for the related discussions.
\end{enumerate}\er


In \cite{HK}, Heinonen and Koskela proved that,  for $Q>1$, quasiconformal mappings between
Ahlfors $Q$-regular metric measure spaces, are quasisymmetric, provided that
the source is a Loewner space and the target space satisfies a quantitative connectivity
condition. See also \cite{hlsw} for a discussion on the quasisymmetry of quasiconformal mappings in $\IR^n$.


\bdefe\label{def-1.3} Let $G\varsubsetneq X$ and $G'\varsubsetneq Y$ be two domains (open and connected), and let $\varphi:[0,\infty)\to [0,\infty)$ be a homeomorphism. We say
that a homeomorphism $f: G\to G'$ is
{\it $\varphi$-semisolid } if $$k_{G'}(f(x),f(y))\leq \varphi(k_G(x,y))$$
for all $x$, $y\in G$ (see Section \ref{sec-2} for the definition of the quasihyperbolic metric).
\edefe

In \cite{GO}, Gehring and Osgood proved that every $K$-quasiconformal mapping $f$ in a domain $G\subsetneq\mathbb{R}^n$ is a $\varphi$-semisolid mapping, where $\varphi(t)=c\max\{t, t^{1/(1-n)}\}$
and $c=c(K, n)$ which means that the constant $c$ depends only on the coefficient $K$ of quasiconformality of $f$ and the dimension $n$ of the Euclidean space $\IR^n$ (see \cite[Theorem 3]{GO}), and thus $f^{-1}$ is also $\varphi$-semisolid since it is $K$-quasiconformal.
In \cite[Theorem 6.12]{TV-2}, Tukia and V\"ais\"al\"a proved that the converse is also true. This implies that
$f$ is $K$-quasiconformal if and only if both $f$ and $f^{-1}$ are $\varphi$-semisolid, where $K$ and $\varphi$ depend on each other and $n$.
This means, by using the terminology in V\"ais\"al\"a's free theory, that quasiconformality and free quasiconformality in $\IR^n$ are quantitatively equivalent.

We introduce the following definition, which is a generalization of convexity:

\bdefe \label{def1.2} For $c\geq 1$, a metric space $X$ is {\it c-quasiconvex} if each pair of points $x$, $y\in X$ can be joined by an curve $\gamma$ with length ${\ell}(\gamma)\leq c|x-y|$.\edefe

We remark that $X$ is convex if and only if $c=1$, and obviously, if $X$ is $c_1$-quasiconvex, then it must be $c_2$-quasiconvex for any $c_2\geq c_1$.


Now, we are ready to state the main result of \cite{HL}, which is the following.

\begin{Thm}\label{Thm-A} $($\cite[Theorem 1.6]{HL}$)$
Let $X$ be a $c$-quasiconvex and complete metric space and $Y$ a
$c'$-quasiconvex metric space. Suppose $f:$ $G\to G'$ is weakly $H$-quasisymmetric, where $G\subsetneq X$ and $G'\subsetneq Y$ are domains. Then there exists a homeomorphism $\varphi:[0,\infty)\to [0,\infty)$ such that $f$ is $\varphi$-semisolid, where $\varphi=\varphi_{c, c', H}$.
\end{Thm}

An open problem given in \cite{HL} can be expressed as follows:
Let $X$ be a $c$-quasiconvex and complete metric space and let $Y$ be a
$c'$-quasiconvex metric space. Suppose that $f:$ $G\to G'$ is homeomorphic, where $G\subsetneq X$ and $G'\subsetneq Y$ are domains
and that there exists a homeomorphism $\varphi:$ $[0,\infty) \to [0,\infty)$ such that $f$ is $\varphi$-semisolid in $G$. Is $f$ weakly $H$-quasisymmetric from $G$ onto $G'$, where $H=H(c, c', \varphi)$?

The following example shows that the answer to this problem is negative.

\beg\label{exa-1} Let $G=\IC\backslash \{0\}$ and $f:$ $G\to G$ be the following inversion:
$$f(z)=\frac{z}{|z|^2}.$$ Then the following statements hold.\ben
\item\label{fri-4}
$f$ is $1$-semisolid in every subdomain of $G$;
\item\label{fri-5}
$f$ is not weakly $H$-quasisymmetric for any $H\geq 1$.
\een
\eeg
Here and thereafter, except when stated otherwise,  $G$ itself is always regarded as a subdomain of $G$. From the proof of Theorem \Ref{Thm-A} in \cite{HL}, we know that Theorem \Ref{Thm-A} also holds in every subdomain of $G$.  Therefore, the problem in \cite{HL} can also be reformulated in the following form.

\bprob\label{prob-1}
Let $X$ be a $c$-quasiconvex and complete metric space and let $Y$ be a
$c'$-quasiconvex metric space. Suppose that $f:$ $G\to G'$ is homeomorphic, where $G\subsetneq X$ and $G'\subsetneq Y$ are domains
and that there exists a homeomorphism $\varphi:$ $[0,\infty) \to [0,\infty)$ such that $f$ is $\varphi$-semisolid in every subdomain of $G$. Is $f$ weakly $H$-quasisymmetric, where $H=H(c, c', \varphi)$?
\eprob

The authors of \cite{HL} studied Problem \ref{prob-1} themselves and obtained the following result.

\begin{Thm}\label{Thm-B}  $($\cite[Theorem 1.8]{HL}$)$
Let $X$ be a $c$-quasiconvex and complete metric space, $G\subsetneq X$ a
non-cut-point domain, and $G'\subsetneq Y$ a domain. Suppose
$f:$ $G\to G'$ is a homeomorphism. If for any subdomain $D\subset G$ and for any $x$ and $y \in D,$
$$k_{f(D)}(f(x), f(y))\leq \varphi(k_D(x, y)),$$
where $\varphi$ is an increasing function, then $f$ is a $K$-quasiconformal
mapping in $G$ with $K=K(c, \varphi)$.
\end{Thm}

Again, it follows from Example \ref{exa-1} that the answer to Problem \ref{prob-1} is negative.
In the following, we shall further develop the ideas of Theorem \Ref{Thm-A} and Problem \ref{prob-1}. In order to state our main result of this paper, two more concepts are required.

\bdefe
Let $G\varsubsetneq X$ and $G'\varsubsetneq Y$ be domains. A homeomorphism $f:$ $G\to G'$ is said to be
\begin{enumerate}\item
{\it $q$-locally $\eta$-quasisymmetric} for some $0<q<1$ if there is a homeomorphism $\eta : [0,\infty) \to [0,\infty)$ such that
$$ |x-a|\leq t|x-b|  \;\;  \mbox{implies} \;\;  |f(x)-f(a)| \leq \eta(t)|f(x)-f(b)|$$
for every $t>0$, for each triplet $x$, $a$, $b$ of points in $\mathbb{B}^{G}(z,q\delta_{G}(z))$  (see Section \ref{sec-3} for the definition) and for all $z\in G$;

\item
{\it $q$-locally weakly $H$-quasisymmetric} for some $0<q<1$ if
$$|x-a|\leq |x-b|\;\;  \mbox{ implies}\;\;   |f(x)-f(a)| \leq H|f(x)-f(b)|$$
for each triplet $x$, $a$, $b$ of points in $\mathbb{B}^G(z,q\delta_{G}(z))$ and for all $z\in G$.\end{enumerate}\edefe

\br\label{rem-1.2}
\begin{enumerate}
\item\label{hwz-1.2}
The weak quasisymmetry implies the  locally weak  quasisymmetry;
\item\label{hwz-12}
Both the locally weak quasisymmetry  and the weak quasisymmetry imply the quasiconformality.
\end{enumerate}
\er


Now, we are ready to state our main result.

\bthm\label{thm-1}
Suppose that $X$ is a $c$-quasiconvex and complete metric space and $Y$ is a $c'$-quasiconvex metric space and that $f:$ $G \to G'$ is homeomorphic, where
$G \varsubsetneq X$ and $G'\varsubsetneq Y$ are domains. Then the following statements are equivalent:
\ben
\item\label{(1)}
there is a homeomorphism $\varphi:$ $[0,\infty) \to [0,\infty)$ such that $f$ is $\varphi$-semisolid in every subdomain of $G$;
\item\label{(2)}
there are constants $\mu>0$ and $0<\alpha<1$ such that $f$ is $\varphi_1$-semisolid in every subdomain of $G$, where $\varphi_1(t)=\mu \max\{t^{\alpha},\;t\}$;
\item\label{(3)}
there are constants $H\geq 1$ and $0<q<1$ such that $f$ is $q$-locally weakly $H$-quasisymmetric in $G$,
\een
where the constants $\mu$, $\alpha$, $H$, $q$ and the control function $\varphi$ depend on each other and the constants $c$ and $c'$.
\ethm

\br\label{rem-1.3} $(i)$ It follows from Remark \ref{rem-1.2} \eqref{hwz-1.2} that the implication from \eqref{(3)} to \eqref{(1)} in Theorem \ref{thm-1} is a generalization of Theorem \Ref{Thm-A}.

$(ii)$ It follows from the discussions in Section \ref{sec-4} that the condition ``$Y$ being $c'$-quasiconvex" is unnecessary in the proof of the implication  from \eqref{(1)} to \eqref{(3)} in Theorem \ref{thm-1}. Hence
Remark \ref{rem-1.2} \eqref{hwz-12} shows that the implication from \eqref{(1)} to \eqref{(3)} in Theorem \ref{thm-1} is a generalization of Theorem \Ref{Thm-B}. Further, this implication implies that the assumption ``$G$ being non-point-cut" in Theorem \Ref{Thm-B} is redundant.

$(iii)$ The equivalence between \eqref{(1)} and \eqref{(3)} in Theorem \ref{thm-1} together with Example \ref{exa-1} can be regarded as a complete answer to Problem \ref{prob-1}.

$(iv)$ As we mentioned in the paragraph next to Definition \ref{def-1.3}, the implication from \eqref{(1)} to \eqref{(2)} in Theorem \ref{thm-1} was first discussed by Gehring and Osgood \cite[Theorem 3]{GO} for the Euclidean
case $X = Y = \IR^n$ with $\mu=\mu(K, n)$ and $\alpha=K^{\frac{1}{1-n}}$ when $f$ is $K$-quasiconformal. Later, Anderson, Vamanamurthy and Vuorinen \cite[Theorem 1.12]{AVV} proved a
version of this result that avoids an explicit reference to the dimension $n$, i.e., $\alpha=K^{\frac{1}{1-n}}$, but $\mu=\mu(K)$, i.e. $\mu$ is independent of $n$. In 1999, V\"ais\"al\"a generalized \cite[Theorem 3]{GO} and \cite[Theorem 1.12]{AVV} to the setting of Banach spaces \cite[Theorem 12.3]{Vai-5}. Obviously, the implication from \eqref{(1)} to \eqref{(2)} in Theorem \ref{thm-1} is a generalization of
\cite[Theorem 12.3]{Vai-5}.
\er

Naturally, the above considerations lead one the following problem:
Under the assumptions of Problem \ref{prob-1}, suppose further that there exists a homeomorphism $\varphi:$ $[0,\infty) \to [0,\infty)$ such that $f$ is $\varphi$-semisolid in $G$.
Are there constants $H\geq 1$ and $0<q<1$ such that $f$ is $q$-locally weakly $H$-quasisymmetric in $G$?

The following example indicates that the answer to this problem is still negative.

\beg\label{exa-2} Let $G=\{z=x+iy\in \IC:\; y>0\}$ and $f:$ $G\to G$ be the following homeomorphism:
$$
f(z)=\begin{cases}
\displaystyle \; x+iy,\;\;\;\;\;\;\;\;\;\;\;\;\;\text{if}\;\; x\leq 0,\\
\displaystyle \; x+i(x+1)y,\;\;\text{if}\;\;  x\geq 0,\; 0<y\leq 1,\\
\displaystyle \; x+i(x+y), \;\;\;\,\text{if}\;\;  x\geq 0,\; y> 1.
\end{cases}$$
 Then the following statements hold:\ben
\item[(a)]
$f$ is $\varphi$-semisolid in G, where $\varphi(t)=\sqrt{3}t$;
\item[(b)]
$f$ is not $\psi$-semisolid in every proper subdomain of $G$ for any homeomorphism $\psi:$ $[0, \infty)\to [0, \infty)$; and
\item[(c)]
$f$ is not $q$-locally weakly $H$-quasisymmetric in $G$ for any $H\geq 1$ and $0<q<1$.
\een
\eeg

We remark that Example \ref{exa-2} was constructed by V\"ais\"al\"a \cite[Example 7.4]{Vai-5} when he showed the conclusion that semisolidity does not imply solidity.\medskip




As a direct consequence of Theorem \ref{thm-1}, we see from \cite[Theorem 3.13]{HK} that the necessity of quasiconformality of homeomorphisms in \cite{S-W} is also sufficient when the space is complete and the image space is unbounded as the following result illustrated.

\bthm
Suppose that $f:$ $G\to G'$ is a homeomorphism between two domains $G\subsetneq X$ and
$G'\subsetneq Y$ of globally $Q$-bounded geometry $(Q > 1)$, that $X$ is complete, and that $Y$ is unbounded.
Then the following statements are equivalent:
\ben
\item
$f$ is $K$-quasiconformal;
\item
There are constants $\mu>0$ and $0<\alpha<1$ such that $f$ is $\varphi_1$-semisolid in every subdomain of $G$, where $\varphi_1(t)=\mu \max\{t^{\alpha},\;t\}$;
\item
There are constants $H\geq 1$ and $0<q<1$ such that $f$ is $q$-locally weakly $H$-quasisymmetric in $G$,
\een
where the constants $K$, $H$, $q$, $\mu$ and $\alpha$ depend on each other and the data of the spaces $X$ and $Y$.
\ethm

It is easy to see that the composition of two $\varphi$-quasisymmetric mappings is also a $\varphi$-quasisymmetric mapping, see \cite[Theorem 2.2]{TV}.
But, in \cite{TV}, Tukia and V\"ais\"al\"a have indicated that the composition of two weakly quasisymmetric mappings need not be weakly quasisymmetric.
So, it is significant to study whether the composition of two locally
weakly quasisymmetric mappings is quasisymmetric or quasiconformal.
As another application of Theorem \ref{thm-1}, we obtain the following.

\begin{thm}\label{thm-3} Suppose that $X_i$ is a $c_i$-quasiconvex and complete metric space $($$i=1, 2$$)$ and that $X_3$ are $c_3$-quasiconvex metric spaces. For domains $G_i\varsubsetneq X_i$ $($$i=1, 2, 3$$)$, if $f:$ $G_1 \to G_2$ is a $q_1$-locally weakly $H_1$-quasisymmetric  mapping and $g:$ $G_2\to G_3$ is a $q_2$-locally weakly $H_2$-quasisymmetric mapping, then\begin{enumerate}
\item the composition $g\circ f$ is $q$-locally $H$-quasisymmetric;
\item the composition $g\circ f$ is also $H$-quasiconformal,
\end{enumerate}
where the constants $q$ and $H$ depend only on $c_1$,
$ c_2$, $c_3$, $q_1$, $q_2$, $H_1$ and $H_2$.
\end{thm}

The rest of this paper is organized as follows. In Section \ref{sec-2}, we shall prove the two examples: Examples \ref{exa-1} and \ref{exa-2}, and some basic results with respect to quasihyperbolic metric and length metric will be shown in Section \ref{sec-3}. The next two sections will be devoted to the proof of the equivalence between \eqref{(1)} and \eqref{(3)} in Theorem \ref{thm-1}. The proof of the equivalence between \eqref{(1)} and \eqref{(2)} in Theorem \ref{thm-1} will be presented in Section \ref{sec-6}, and in Section \ref{sec-7}, we shall prove Theorem \ref{thm-3}.

\section{Properties of the examples}\label{sec-2}

The purpose of this section is to prove the two examples constructed in the first section.

\subsection{The proof of Example \ref{exa-1}}
It follows from \cite[Theorem 5.11]{Vai-5} that $f$ is an isometry with respect to the quasihyperbolic metric. Hence the first statement holds.

For the proof of \eqref{fri-5}, we let $a=\frac{1}{t}$, $x=1$ and $b=t$, where $t>1$. Suppose on the contrary that $f$ is weakly $H$-quasisymmetric for some $H>0$. Since $|a-x|\leq |b-x|$, we have $$\frac{|f(a)-f(x)|}{|f(b)-f(x)|}\leq H,$$
which implies that for any $t>1$,
$$t\leq H.$$ This is impossible.
Hence the proof is complete.\qed

\subsection{The proof of Example \ref{exa-2}}

Obviously, the first statement $(a)$ follows from \cite[Example 7.4]{Vai-5}.

Now, we prove $(c)$.
Suppose on the contrary that $f$ is $q$-locally weakly $H$-quasisymmetric in $G$ for some $H\geq 1$ and $0<q<1$.
Let $0<\varepsilon<\frac{1}{2}$, $O=(n, \frac{1}{2})$ and $\mathbb{B}(O, \varepsilon)=\{z\in \IC:\; |z-O|<\epsilon\}$. Then $\mathbb{B}(O, \varepsilon)\subset G$, and so
 $$q\mathbb{B}(O, \varepsilon)=\mathbb{B}(O, q\varepsilon)\subset G.$$

Let
$$a=(n, \frac{1}{2}+\frac{1}{2}q\varepsilon)\;\;\mbox{and}\;\; b=(n-\frac{1}{2}q\varepsilon, \frac{1}{2}).$$
Then we easily know that
$$a, O, b\in q\mathbb{B}(O, \varepsilon)\;\;\mbox{and}\;\; |a-O|=|b-O|.$$

But elementary computations show that
$$\frac{|f(a)-f(O)|}{|f(b)-f(O)|}=\frac{2\sqrt{5}}{5}(n+1),$$ which implies that $f$ is not weakly $H$-quasisymmetric in $q\mathbb{B}(O, \varepsilon)$ for any $H\geq 1$.

Next, we show $(b)$. Suppose that $f$ is $\psi$-semisolid in every proper subdomain of $G$ for some homeomorphism $\psi:$ $[0, \infty)\to [0, \infty)$. Then $(a)$ and
Theorem \ref{thm-1} or \cite[Theorem $7.9$]{Vai-5} guarantee that $f$ is $q$-locally weakly $H$-quasisymmetric in $G$ for some $H\geq 1$ and $0<q<1$.
This contradicts $(c)$. Hence the proof of the example is complete.\qed

\section{Auxiliary results}\label{sec-3}

\subsection{Basic geometric results in metric spaces}

We always denote the open metric ball with center $x\in X$ and radius $r>0$ by
$$\mathbb{B}(x,r)=\{z\in X:\; |z-x|<r\}$$
and for $\lambda>0$, $$\lambda \mathbb{B}(x,r)=\{ z\in X:\; |z-x|<\lambda r\}.$$

For a set $A$ in $X$, we always use $\partial A$ (resp. $\overline{A}$) to denote the boundary (resp. the closure) of $A$.

By a curve, we mean any continuous function $\gamma:$ $[a,b]\to X$. The length of $\gamma$ is denoted by
$$\ell(\gamma)=\sup\Big\{\sum_{i=1}^{n}|\gamma(t_i)-\gamma(t_{i-1})|\Big\},$$
where the supremum is taken over all partitions $a=t_0<t_1<t_2\ldots<t_n=b$. The curve is {\it rectifiable} if $\ell(\gamma)<\infty$.
In particular, if the metric is taken to be the quasihyperbolic metric, the length of $\gamma$ is denoted by $\ell_{qh}(\gamma)$.

The length function associated with a rectifiable curve $\gamma$: $[a,b]\to X$ is $s_{\gamma}$: $[a,b]\to [0, \ell(\gamma)]$, given by
$s_{\gamma}(t)=\ell(\gamma|_{[a,t]})$.
For any rectifiable curve $\gamma:$ $[a,b]\to X$, there is a unique curve $\gamma_s:$ $[0, \ell(\gamma)]\to X$ such that $\gamma=\gamma_s\circ s_{\gamma}$. Obviously,
$\ell(\gamma_s|_{[0,t]})=t$ for $t\in [0, \ell(\gamma)]$. The curve $\gamma_s$ is called the arclength parametrization of $\gamma$.

For a rectifiable curve $\gamma$ in $X$, the line integral over $\gamma$ of each Borel function $\varrho:$ $X\to [0, \infty)$ is
$$\int_{\gamma}\varrho ds=\int_{0}^{\ell(\gamma)}\varrho\circ \gamma_s(t) dt.$$

Let $X$ be a connected metric space and $\Omega$ a nonempty proper subset of $X$. It follows from \cite[Remark 2.2]{HL} that if $\Omega$ is open, then $\partial \Omega\not=\emptyset$.
Further, we have
\blem\label{mon-4}
Suppose $X$ is a connected metric space. For any nonempty set $\Omega\subsetneq X$, $\partial \Omega\not=\emptyset$.
\elem
\bpf
Suppose on the contrary that $\partial \Omega=\emptyset$. Let $u\in \Omega$. Then there is an open set $V\subset X$ such that $u\in V$ and either $V\subset \Omega$ or  $V\subset X\backslash \Omega$. Since $u\in V\cap \Omega$, we see that $V\subset \Omega$. This implies that $\Omega$ is open. By \cite[Remark 2.2]{HL}, it is impossible. Hence the lemma holds.\epf

For $z\in \Omega$, we always use $\delta_{\Omega}(z)$ to denote the distance from $z$ to $\partial \Omega$. We make the convention that $\delta_{\Omega}(z)=0$ if and only if $z\in \partial \Omega$. Then we easily have the following useful inequality: For any $x$, $y\in \Omega$,
$$\delta_{\Omega}(x)\leq \delta_{\Omega}(y)+|x-y|.$$

For $z\in \Omega$, assume that $\delta_{\Omega}(z)>0$. Let $r\in (0, \delta_{\Omega}(z))$. Then the metric ball $\mathbb{B}(z, r)$ is not necessarily contained in $\Omega$, and also the intersection $\mathbb{B}(z, r)\cap \Omega$ is not always connected even when $\Omega$ is a domain. So we always consider the component of the intersection $\mathbb{B}(z, r)\cap \Omega$ containing the center $z$, which is denoted by $\mathbb{B}^{\Omega}(z,r)$, when $\Omega$ is open. Similarly, we use
$\overline{\mathbb{B}}^{\Omega}(z,r)$ (resp. $\partial \mathbb{B}^{\Omega}(z,r)$) to denote the closure (resp. the boundary) of the component of $\mathbb{B}(z, r)\cap \Omega$ containing the center $z$.
In particular, for $\lambda >0$, we always use
$\lambda \mathbb{B}^{\Omega}(z,r)$ to denote the component $\mathbb{B}^{\Omega}(z, \lambda r)$ provided $\lambda r<\delta_{\Omega}(z)$. Then we have the following.

\begin{lem}\label{lem-3.1} Suppose $X$ is a $c$-quasiconvex metric space and $\Omega \subsetneq X$ is open. For any rectifiably connected set
$D\subset \mathbb{B}(z,r)$ with $z\in D\cap \Omega$, if $r\leq \delta_{\Omega}(z)$, then $$D \subset \mathbb{B}^{\Omega}(z,r) \subset \Omega.$$
\end{lem}
\bpf Obviously, it suffices to prove that $D\subset \Omega$. We show this by contradiction.
Suppose on the contrary that $D$ is not contained in $\Omega$. Then there exists a point $u\in D\subset \mathbb{B}(z,r)$, but $u \not\in \Omega$.
Let $\gamma$ be a rectifiable curve in $D$ joining the points $z$ and $u$. Since $\partial \Omega\not=\emptyset$ and
$u \not\in \Omega$, we see that there must exist a point $u_0$ such that $u_0 \in \gamma \cap \partial \Omega$. Hence
$$|z-u_0|\geq \delta_{\Omega}(z)\geq r> |z-u_0|,$$
which is the desired contradiction.
\epf

For two domains $D$ and $G$ in a space, in general, the assuption $D\subset G$ implies $\mathbb{B}^D(z,\delta_D(z))\subset \mathbb{B}^{\Omega}(z,\delta_{\Omega}(z))$. This basic and very useful property holds for many types of spaces, for example in Banach spaces. But this property is not valid in metric spaces as the following example shows.
\beg\label{mon-3}
Let
$$X=l_1\cup l_2\cup l_3\cup l_4, \;\;\Omega=X\backslash l_5\;\;\mbox{and}\;\; D=l_1,$$
where
$l_1=\{(x,0): \; -2< x< 2\}$, $l_2=\{(x,1): \; -2< x< 2\}$, $l_3=\{(-2,y): \; 0\leq y\leq 1\}$, $l_4=\{(2,y): \; 0\leq y\leq 1\}$
and $l_5=\{(x,1): \; -1\leq x\leq 1\}$. Then, obviously, we have
\ben
\item
$X$ is $5$-quasiconvex with the Euclidean metric;
\item\label{thu-5}
$D\subsetneq \Omega$, $\delta_D(z)=2$ and $\delta_{\Omega}(z)=\sqrt{2}$, where $z=(0,0)$, and
\item\label{thu-6}
$\mathbb{B}^{\Omega}(z,\delta_{\Omega}(z))\subsetneq \mathbb{B}^D(z,\delta_D(z))$.
\een
\eeg

Our next lemma illustrates that with some constraint, this phenomenon can be avoided.
\begin{lem}\label{thu-3} Suppose $X$ is a $c$-quasiconvex metric space and $\Omega \subsetneq X$ is open.
\ben
\item\label{mon-01}
 For any nonempty subset $P$ in $\Omega$, we have $\delta_P(z)\leq c\delta_{\Omega}(z)$ for any $z\in P$;
\item\label{mon-02}
 For any open subset $D$ in $\Omega$, we have $\mathbb{B}^D\Big(z,\frac{r}{c}\delta_D(z)\Big)\subset \mathbb{B}^{\Omega}(z,r\delta_{\Omega}(z))$ for any $z\in D$, where $0<r\leq 1$.
\een\elem
\bpf
Obviously, the second statement easily follows from the first one by taking $P=D$. So, to prove this lemma, it suffices to show the first statement. Without loss of generality, we may assume that $P\subsetneq \Omega$ and $\delta_P(z)>0$. Let $a\in \partial \Omega$. Then we claim that
\be\label{thu-4}\delta_P(z)\leq c|z-a|.\ee

It follows from Lemma \ref{mon-3} that $\partial P\not=\emptyset$. To prove this inequality, we consider two cases.
 For the first case where $a\in \partial P$, obviously, $$\delta_P(z)\leq |z-a|.$$
For the remaining case, that is, $a\notin \partial P$, we know that $a\in X\backslash \overline{P}$. Let $\alpha$ denote a curve in $X$ connecting $a$ and $z$ such that
$$\ell(\alpha)\leq c|z-a|.$$

Since $z\in P$, obviously, $\alpha\cap \partial P\not=\emptyset$. Assume that $b\in \alpha\cap \partial P$. Then
$$\delta_P(z)\leq |z-b| <\ell(\alpha)\leq c|z-a|.$$ Hence the inequality \eqref{thu-4} holds.\medskip

By taking the infimum over all $a\in \partial \Omega$ in \eqref{thu-4}, we easily know that the statement \eqref{mon-01} in the lemma is true, and so the proof of the lemma is complete.
\epf

Further, we have the following result.

\begin{lem}\label{lem-3.2} Suppose $X$ is a $c$-quasiconvex metric space and $\Omega \subsetneq X$ is open.
\begin{enumerate}
\item\label{hwz-130} Suppose $z\in \Omega$ and $0<r\leq \frac{2}{2+c}\delta_{\Omega}(z)$. Then
$\mathbb{B}(z, r)\subset \Omega$;
\item\label{hwz-131} Suppose $z\in \Omega$ and $0<r\leq \frac{1}{1+c}\delta_{\Omega}(z)$. Then
for any $x$ and $y\in \mathbb{B}(z, r)$, there must exist a curve $\gamma\subset \mathbb{B}^{\Omega}(z, (1+c)r)$ such that
$$\ell(\gamma)\leq c|x-y|;$$
\item\label{hwz-132}
Suppose that $z\in \Omega$ and $0<r\leq \frac{1}{1+c}\delta_{\Omega}(z)$. Then
$$\mathbb{B}(z, r)\subset \mathbb{B}^{\Omega}(z, (1+c)r).$$
In particular, for any $z\in \Omega$ and $0<r\leq \frac{1}{1+c}\delta_{\Omega}(z)$, $$\mathbb{B}(z, r)\subset \mathbb{B}^{\Omega}(z, \delta_{\Omega}(z)).$$
\end{enumerate}
\end{lem}
\bpf
We prove the first statement by contradiction. Suppose there exists some $x\in \mathbb{B}(z, r)$ such that $x\notin \Omega$. By the assumption, we see that there is a curve $\beta\subset X$ connecting $z$ and $x$ such that
$$\ell(\beta)\leq c|x-z|.$$
Then $\ell(\beta)<cr$. We claim that for $w\in \beta$,
\be\label{hwz-110-2}|w-z|  <  \frac{2+c}{2}r.\ee
To show this claim, let $w_0$ be the point in $\beta$ such that $$\ell(\beta_{[z, w_0]})=\ell(\beta_{[w_0, x]}),$$ where $\beta_{[z, w_0]}$ denotes the part of $\beta$ with the endpoints $z$ and $w_0$. According to the position of $w$ in $\beta$, we consider two possibilities. For the first possibility where $w\in \beta_{[z, w_0]}$, we easily see that
$$|w-z|\leq \frac{1}{2}\ell(\beta)< \frac{c}{2}r.$$
For the other possibility, that is, $w\in \beta_{[x, w_0]}$, we can easily get that
$$|w-z|\leq |w-x|+|x-z|< \frac{1}{2}\ell(\beta)+r\leq \frac{2+c}{2}r.$$
Hence \eqref{hwz-110-2} holds, and the claim is proved.\medskip

Now, we proceed with the proof based on \eqref{hwz-110-2}. Since $z\in \Omega$ and $\partial\Omega\not=\emptyset$, we know that $\beta\cap \partial \Omega\not=\emptyset$. Assume that $z_0\in \beta\cap \partial \Omega$. Then we deduce from \eqref{hwz-110-2} that
$$|z_0-z|\geq \delta_{\Omega}(z)\geq \frac{2+c}{2}r>|z_0-z|.$$ This is the desired contradiction, which shows that the statement \eqref{hwz-130} in the lemma is true.\medskip

To get proofs of the second and the third statements in the lemma, for any $x$ and $y\in \mathbb{B}(z, r)$, we let $\gamma$ denote a curve in $X$ joining $x$ to $y$ such that $${\ell}(\gamma)\leq c|x-y|.$$
Then ${\ell}(\gamma)<2cr.$ Now, we check that
\be\label{hwz-140} \gamma\subset \mathbb{B}(z, (1+c)r)\cap \Omega.\ee
Since for all $w\in \gamma$,
\be\label{hwz-110} |w-z| \leq  \min\{|w-x|, |w-y|\}+r
\leq \frac{1}{2}{\ell}(\gamma) + r<  (1+c)r,\ee which implies that $\gamma\subset \mathbb{B}(z,(1+c)r)$.

It remains to prove the inclusion $\gamma\subset \Omega$. Again, we show this by contradiction. Suppose that $\gamma$ is not contained in $\Omega$. Since it follows from the statement \eqref{hwz-130} that $x$ and $y\in \Omega$, there must be a point $u_0$ which is contained in the intersection $\gamma\cap \partial \Omega$, and so \eqref{hwz-110} implies
$$\delta_{\Omega}(z)\leq |u_0-z| <  (1+c)r\leq \delta_{\Omega}(z).$$ This is the desired contradiction. Hence $\gamma\subset \mathbb{B}(z, (1+c)r)\cap \Omega$.\medskip

For any point $v\in \gamma$, let $\zeta$ denote a curve in $X$ connecting $z$ and $v$ such that $${\ell}(\zeta)\leq c|v-z|.$$
By replacing $x$ (resp. $y$) by $z$ (resp. $v$) in \eqref{hwz-140}, we can know that
$$\zeta\subset \mathbb{B}(z, (1+c)r)\cap \Omega.$$ The arbitrariness of $v\in \gamma$ makes sure that $$\gamma\subset \mathbb{B}^{\Omega}(z, (1+c)r).$$ So the second statement in the lemma holds.

Obviously, the statement \eqref{hwz-132} follows from \eqref{hwz-140} and the arbitrariness of the points $x$ and $y$ in $\mathbb{B}(z, r)$, and hence the proof of our lemma is complete.
\epf

In the following, we get several applications of Lemma \ref{lem-3.2}, which will be useful in the discussions later on.

\begin{lem}\label{lem-3.3} Suppose $X$ is a $c$-quasiconvex metric space and $\Omega \subsetneq X$ is open. For any metric ball $\mathbb{B}(z, r)$ with $z\in \Omega$ and $0<r\leq \delta_{\Omega}(z)$, let $B=\mathbb{B}^{\Omega}(z, r)$. Then \ben
\item\label{thu-7}
$\delta_B(z)=r;$
\item\label{thu-8}
for any $u\in \overline{B}$, $|u-z|=r$ if and only if $u\in \partial B$.
\een
\end{lem}
\bpf
By Lemma \ref{mon-3}, we know that $\partial B\not=\emptyset$. Obviously, we only need to prove \eqref{thu-7} since \eqref{thu-8} is a direct consequence of \eqref{thu-7}. We prove \eqref{thu-7} by contradiction.

Suppose $\delta_B(z)<r$.  Then there must exist $z_0\in \partial B$ such that
$$|z_0-z|\leq \delta_B(z)+\frac{r-\delta_B(z)}{2}<r.
$$
Since $z_0\in \partial B$ and $B\subset \Omega$, we know that $z_0\in \overline{\Omega}$. Obviously, $z_0\notin \partial \Omega$, because otherwise,
$\delta_{\Omega}(z)\leq |z-z_0|<r$. Hence $z_0\in \Omega$.

Let $r_1=r-|z_0-z|$. Then $r_1>0$ and
$$\delta_{\Omega}(z_0)>\delta_{\Omega}(z)-|z-z_0|\geq r_1,$$
and so it follows that $$\mathbb{B}(z_0, r_1)\subset \mathbb{B}(z, r).$$

Further, let $r_2=\frac{1}{1+c}r_1$. Then Lemma \ref{lem-3.2} \eqref{hwz-132} implies
$$\mathbb{B}(z_0, r_2)\subset \mathbb{B}^{\Omega}(z_0, r_1).$$

Since $z_0\in \partial B$, we see that there is a sequence $\{x_i\}_{i=1}^{\infty}\subset B$ such that $\lim_{i\to \infty}x_i=z_0$. Then for all sufficiently large $i$, $x_i\in \mathbb{B}(z_0, r_2)$, which implies $\mathbb{B}^{\Omega}(z_0, r_1)\cap B\not=\emptyset$. Hence it follows from the fact $\mathbb{B}(z_0, r_1)\subset \mathbb{B}(z, r)$ that $\mathbb{B}^{\Omega}(z_0, r_1)\subset B,$ and so
$$\mathbb{B}(z_0, r_2)\subset B.$$
This is impossible since $z_0\in \partial B$. Hence \eqref{thu-7} holds, and so the proof of lemma is finished.
\epf

\blem\label{tue-6}
Suppose $X$ is a $c$-quasiconvex metric space and $\Omega \subsetneq X$ is open. Then for any $0<\mu \leq 1$, $\mathbb{B}^{\Omega}(z, \mu\delta_{\Omega}(z))$ is a domain in $\Omega$ for any $z\in \Omega$.
\elem
\bpf
Obviously, it suffices to show the openness of $\mathbb{B}^{\Omega}(z, \mu\delta_{\Omega}(z))$. Since for any $w\in B=\mathbb{B}^{\Omega}(z, \mu\delta_{\Omega}(z)),$
$w\notin \partial B$. Otherwise,
Lemma \ref{lem-3.3} guarantees that
$$\mu\delta_{\Omega}(z)=\delta_B(z)\leq |z-w|<\mu\delta_{\Omega}(z).$$ Hence this lemma is proved.
\epf

By Lemma \ref{tue-6}, we have the following.

\blem\label{tue-7}
Suppose $X$ is a $c$-quasiconvex metric space. Then $X$ must be locally connected.
\elem
\bpf
Let $z$ be any point in $X$, and let $U$ be any open set in $X$ with $z\in U$. To prove this lemma, we only need to find a domain $W$ in $U$ with $z\in W$.

Since $U$ is open and $X$ is $c$-quasiconvex, it follows from \cite[page 120]{JM} that there is $r>0$ such that $\mathbb{B}(z, r)\subsetneq U$. Let
$$V=\mathbb{B}(z,r)\;\;\mbox{and}\;\;W=\mathbb{B}^V(z,\mu\delta_V(z)),$$
 where $\mu\delta_V(z)<r$.
Then we deduce from Lemma \ref{tue-6} that this $W$ is our needed. Hence Lemma \ref{tue-7} is true.
\epf

\subsection{Quasihyperbolic metric}

The {\it quasihyperbolic length} of a rectifiable curve or a path
$\gamma$ in the metric in a domain $G \varsubsetneq X$ is the number:

$$\ell_{k_G}(\gamma)=\int_{\gamma}\frac{|dz|}{\delta_{G}(z)}.
$$
 For any $z_1$, $z_2$ in $G$, the {\it quasihyperbolic distance}
$k_G(z_1,z_2)$ between $z_1$ and $z_2$ is defined by
$$k_G(z_1,z_2)=\inf\{\ell_{k_G}(\gamma)\},
$$
where the infimum is taken over all rectifiable curves $\gamma$
joining $z_1$ to $z_2$ in $G$.

Gehring and Palka \cite{GP} introduced the quasihyperbolic metric of
a domain in $\IR^n$. For the basic properties of this metric we refer to \cite{GO}. Recall that a curve $\gamma$ from $z_1$ to
$z_2$ is a {\it quasihyperbolic geodesic} if
$\ell_{k_G}(\gamma)=k_G(z_1,z_2)$. Each subcurve of a quasihyperbolic
geodesic is obviously a quasihyperbolic geodesic. It is known that a
quasihyperbolic geodesic between any two points in $E$ exists if the
dimension of $X$ is finite, see \cite[Lemma 1]{GO}. This is not
true in arbitrary metric spaces (cf. \cite[Example 2.9]{Vai4}).

We establish comparison results between the metrics $|\cdot|$ and $k_G$ in a $c$-quasiconvex metric space. It is a modified version of \cite[Theorems 2.7 and 2.8]{HL}.


\begin{lem}\label{lem-3.4} Let $X$ be a {\it c-quasiconvex} metric space and let $G \varsubsetneq X$ be a domain.

\noindent
$(1)$ For all $x$, $y\in G$, \be\label{sun-1} |x-y|\leq (e^{k_G(x,y)}-1)\delta_{G}(x);\ee

\noindent
$(2)$ Suppose $z\in G$ and $0<t< 1$. Then for  $x,$ $y\in \overline{\mathbb{B}}^{G}\big(z, \frac{t}{2c}\delta_{G}(z)\big)$,
\be\label{vvm-1} \frac{c}{c+t}\cdot\frac{|x-y|}{\delta_{G}(z)}\leq k_G(x,y) \leq\frac{c}{1-\frac{1+c}{2c}t}\cdot\frac{|x-y|}{\delta_{G}(z)};\ee

\noindent
$(3)$ Suppose that $x$, $y\in G$ and either $|x-y|\leq \frac{1}{3c}\delta_{G}(x)$ or $k_G(x,y)\leq 1$. Then
\be\label{vvm-2} \frac{1}{2}\frac{|x-y|}{\delta_{G}(x)}< k_G(x,y) \leq 3c\frac{|x-y|}{\delta_{G}(x)}.\ee
\end{lem}

\bpf
The inequality \eqref{sun-1} follows from \cite[Theorem 2.7]{HL}. In the following, we prove \eqref{vvm-1}  and \eqref{vvm-2}.

Since $X$ is a $c$-quasiconvex metric space, we see from Lemma \ref{lem-3.2} \eqref{hwz-131} that there is a curve $\gamma$ in $\mathbb{B}^{G}(z, \frac{1+c}{2c}t\delta_{G}(z))$ joining $x$ and $y$ with $${\ell}(\gamma)\leq c|x-y|.$$
Since for $w\in \gamma$,
$$\delta_{G}(w)\geq \delta_{G}(z)-|w-z|\geq \Big(1-\frac{1+c}{2c}t\Big)\delta_{G}(z),$$
we see that
$$k_G(x,y)\leq \int_{\gamma}\frac{|dw|}{\delta_{G}(w)}\leq \frac{c}{1-\frac{1+c}{2c}t}\cdot\frac{|x-y|}{\delta_{G}(z)}.$$
This shows that the right-side inequality of \eqref{vvm-1} is true. Now, we give a proof of the left-side inequality of \eqref{vvm-1}.
It follows from \cite[Observation 2.6]{HL} that $G$ is rectifiably connected.
Hence it suffices to prove that for any rectifiable curve $\alpha\subset G$ joining $x$ and $y$,
\beq\label{sun-2}\nonumber
\ell_{k_G}(\alpha)\geq \frac{c}{c+t}\cdot\frac{|x-y|}{\delta_{G}(z)}.
\eeq
To this end, let $\alpha$ be a rectifiable curve in $G$ joining $x$ and $y$.
We divide the proof into two cases. For the first case where $\alpha\subset \overline{\mathbb{B}}^{G}(z,\frac{t}{c}\delta_{G}(z))$, we see that
\be\label{japan-34}\delta_{G}(w)\leq |w-z|+\delta_{G}(z)\leq \Big(1+\frac{t}{c}\Big)\delta_{G}(z)\ee for all $w\in \alpha$,
and hence
$$\ell_{k_G}(\alpha)=\int_{\alpha}\frac{|dw|}{\delta_{G}(w)}\geq \frac{c}{c+t}\cdot\frac{|x-y|}{\delta_{G}(z)},$$ as required.

For the remaining case, that is, $\alpha\nsubseteq \overline{\mathbb{B}}^{G}(z,\frac{t}{c}\delta_{G}(z))$, obviously, $\alpha$ has two sub-curves $\alpha_{1}$ and $\alpha_{2}$
 in $\overline{\mathbb{B}}^{G}(z, \frac{t}{c}\delta_{G}(z))$ joining the sets $\partial\mathbb{B}^{G}(z, \frac{t}{2c}\delta_{G}(z))$ and $\partial\mathbb{B}^{G}(z, \frac{t}{c}\delta_{G}(z))$.
Again, \eqref{japan-34} implies
$$\ell_{k_G}(\alpha)\geq \int_{\alpha_{1}\cup\alpha_{2}}\frac{|dw|}{\delta_{G}(w)}\geq \frac{c}{c+t}\cdot\frac{|x-y|}{\delta_{G}(z)},$$
since $\ell(\alpha_i)\geq \frac{t}{2c}\delta_{G}(z)\geq \frac{|x-y|}{2}$ for $i=1$, $2$.
Hence the proof of \eqref{vvm-1} is complete.
\medskip

To prove \eqref{vvm-2}, we first consider the case where $|x-y|\leq \frac{1}{3c}\delta_{G}(x)$. By taking  $t= \frac{2}{3}$ in \eqref{vvm-1}, obviously, we have that
 $$\frac{|x-y|}{2\delta_{G}(x)}< \frac{c|x-y|}{(c+\frac{2}{3})\delta_{G}(x)}\leq k_G(x,y)
\leq \frac{3c^{2}}{2c-1}\cdot\frac{|x-y|}{\delta_{G}(x)}\leq 3c\frac{|x-y|}{\delta_{G}(x)}.$$
For the remaining case, that is, $|x-y|> \frac{1}{3c}\delta_{G}(z)$ and $k_G(x,y)\leq 1$, the right-side inequality of \eqref{vvm-2} is obvious. Since for $r\in (0,1]$, $$e^r-1< 2r,$$ it follows from \eqref{sun-1} that $$|x-y|\leq (e^{k_G(x,y)}-1)\delta_{G}(x)< 2k_G(x,y)\delta_{G}(x),$$ which implies the left-side inequality of \eqref{vvm-2} is true too. Hence the proof of \eqref{vvm-2} is  complete.
\epf

The following result is on the quasiconvexity of a domain $G$ in $X$ with respect to the quasihyperbolic metric $k_G$.

\begin{lem}\label{lem-3.5}  Let $X$ be a $c$-quasiconvex metric space and let $G \varsubsetneq X$ be a domain.
\ben
\item\label{sun-3}Let $\gamma$ be a rectifiable path in $G$.
Then  $\ell_{k_G}(\gamma)$ is the length of $\gamma$ in the metric space $(G,k_G)$, i.e.,
$\ell_{k_G}(\gamma)=\ell_{qh}(\gamma)$;
\item\label{sun-4}
Then the space $(G,k_G)$ is  $\lambda$-quasiconvex for all $\lambda >1.$\een
\end{lem}
\bpf The proof of the statement \eqref{sun-3} easily follows from \eqref{vvm-1} in Lemma \ref{lem-3.4} and
\cite[Lemma 2.6]{BHK}, and the second statement easily follows from the first one.
\epf

We remark that in Lemma \ref{lem-3.5} and in what follows,  the notation $(G, e)$ is applied to emphasize the related metric $e$.

\subsection{Length metric}

In the following, we always use $d$ to denote the length metric in $X$ with respect to $|\cdot|$ and $k'_G$ to denote the quasihyperbolic metric with respect to $d$ in $G$, where $G\subsetneq X$ is a domain. Then we have the following analog of statements $(i)$ and $(ii)$ in \cite[Lemma 9.2]{HSX} in the setting of $c$-quasiconvex metric spaces, i.e. without the assumption that $X$ is
proper.

\blem\label{lem-3.6} Suppose $X$ is a $c$-quasiconvex metric space and $G \subsetneq X$ is a domain. Then
\begin{enumerate}
\item\label{wed-3}
for all $x,$ $y\in X$,
$|x-y|\leq d(x, y) \leq c|x-y|;$ and
\item\label{wed-4}
for all $x,$ $y\in G$,
$\frac{1}{c}k_G(x, y)\leq k'_G(x, y) \leq c k_G(x, y).$
\end{enumerate}
\elem
\bpf
To prove this lemma, it suffices to verify the inequality \eqref{wed-4} in the lemma since the proof of \eqref{wed-3} is obvious.

By the inequality \eqref{wed-3} in the lemma, it follows that for any curve $\gamma\subset X$, $$\ell(\gamma)\leq \ell_d(\gamma)\leq c\ell(\gamma),$$
 where $\ell_d(\gamma)$ denotes the arclength of $\gamma$ with respect to $d$. Let $\delta'_G(x)$ denote the distance from $x$ to the boundary of $G$ with respect to $d$.
Then we have that for all $x\in G$,
 \be\label{wed-5}  \delta_{G}(x)\leq \delta'_G(x)\leq c\delta_{G}(x).\ee

For the proof, we let $\varepsilon$ be a positive number. Then there are points $a,$ $b\in \partial G$ such that
$$|x-a|\leq \delta_{G}(x)+\frac{\varepsilon}{c}\;\;\mbox{and}\;\; d(x,b)\leq \delta'_G(x)+\varepsilon,$$ respectively.

By \eqref{wed-3}, we easily know that
$$\delta_{G}(x)\leq |x-b|\leq d(x,b)\leq \delta'_G(x)+\varepsilon$$
and $$\delta'_G(x)\leq d(x,a)\leq c|x-a|\leq c\delta_{G}(x)+\varepsilon.$$
Thus the inequality \eqref{wed-5} easily follows from the arbitrariness of $\varepsilon$.\medskip

Now, we can finish the proof by using \eqref{wed-5}. Let $|d_d(z)|$ denote the arclength differential with respect to the length metric $d$. Then
$$|dz|\leq |d_d(z)|\leq c|dz|.$$

Since for all $x,$ $y\in G$ and for any $h>0$, there is a curve $\gamma$ (resp. $\gamma'$) connecting $x$ and $y$ in $G$ such that $$\ell_{k_{G}}(\gamma)\leq k_G(x,y)+\frac{h}{c}\;\;(\mbox{resp.}\;\; \ell_{k'_G}(\gamma')\leq k'_G(x,y)+\frac{h}{c}),$$ we know that
$$k'_G(x,y)\leq \int_{\gamma}\frac{|d_d(z)|}{\delta'_G(z)}\leq \int_{\gamma}\frac{c|dz|}{\delta_{G}(z)}=c\ell_{k_{G}}(\gamma) \leq ck_G(x,y)+h$$
and
$$k_G(x,y)\leq \int_{\gamma'}\frac{|dz|}{\delta_{G}(z)}\leq \int_{\gamma'}\frac{c|d_d(z)|}{\delta'_G(z)}=c\ell_{k'_G}(\gamma') \leq ck'_G(x,y)+h.$$
Obviously, letting $h\to 0$ yields the inequality \eqref{wed-4} in the lemma, and so the proof of the lemma is complete.\medskip
\epf

\section{The proof of the implication from \eqref{(1)} to \eqref{(3)} in Theorem \ref{thm-1}}\label{sec-4}

This section consists of two subsections. In the first subsection, we shall show the implication from relativity and locally weak quasisymmetry with a constraint on the coefficient of quasiconvexity. In the second subsection, the proof of the implication from \eqref{(1)} to \eqref{(3)} in Theorem \ref{thm-1} will be proved.

\subsection{Relativity and locally weak quasisymmetry}

We start this subsection with the definition of relativity of homeomorphisms.

\bdefe Suppose $f:$ $G\to G'$ is a homeomorphism, where $G\subsetneq X$ and $G'\subsetneq Y$ are domains. Then $f$ is said to be \begin{enumerate}
\item
$(\theta ; t_{0})$-relative if there are $t_0\in (0, 1]$ and a homeomorphism $\theta$: $[0, t_0)\to [0, \infty)$ such that $$\frac{|f(x)-f(y)|}{\delta_{G'}(f(x))} \leq\theta\Big(\frac{|x-y|}{\delta_G(x)}\Big)$$
whenever $x$, $y\in G$ with $|x-y|<t_{0}\delta_G(x)$;
\item
$\theta$-relative if $t_0=1$.
\end{enumerate}\edefe

The following is the main result in this subsection.

\blem\label{lem-4.1}
Suppose that $X$ is a $c$-quasiconvex metric space with $1\leq c\leq \frac{1+\sqrt{3}}{2}$ and that both $G\subsetneq X$ and $G'\subsetneq Y$ are domains. If $f:$ $G\to G'$ is $\theta$-relative in every subdomain of $G$, then $f$ is $q'$-locally $\eta'$-quasisymmetric in $G$ with $q'=q'(c)=\frac{1}{(2+c)^3}$ and $\eta'=\eta'_{\theta,c}$.
\elem
\bpf
Since $\theta:$ $[0, 1)\to [0, \infty)$ is homeomorphic, without loss of generality, we assume that $$M=\theta\Big(\frac{2c}{1+2c}\Big)>1.$$

For distinct points $x$, $a$, $b$ in $\mathbb{B}^G(z,q\delta_{G}(z))$ with $|x-a|=t|x-b|$, where $z\in G$ and $t>0$, to prove that $f$ is locally quasisymmetric in $G$, we need a relationship between $|f(a)-f(x)|$ and $|f(x)-f(b)|.$ To this end, we divide the discussions into three cases according to the location of the parameter $t$.
We first discuss the case where $0<t\leq \frac{2c}{1+2c}$. In this case, we have
\medskip
\bcl\label{tue-1}
Suppose $0<t\leq \frac{2c}{1+2c}$ and $0<q_1\leq \frac{1}{2+c}$. Then \ben
\item\label{wed-6}
$x$ and $a$ must lie in a component of $G\setminus \{b\}$, denoted by $G_{abx}$. Further, $G_{abx}$ is open;
\item\label{wed-7}
$|f(a)-f(x)|\leq \theta(t)|f(x)-f(b)|.$
\een\ecl

We prove \eqref{wed-6} in the claim by contradiction. Suppose on the contrary that $x$ and $a$ lie in two different components of  $G\setminus \{b\}$. By Lemma \ref{lem-3.2} (\ref{hwz-131}), there exists a curve $\alpha \subset G$ joining $x$ and $a$ with
$$\ell(\alpha)\leq c|x-a|.$$
Then from the contrary assumption we know that
$b\in \alpha.$
Moreover, $$|x-b|<\ell(\alpha)\leq c|x-a|=ct|x-b|,$$
and so $$\frac{1}{c}< t \leq \frac{2c}{1+2c}.$$ This is the desired contradiction since we have assumed that $c\leq \frac{\sqrt{3}+1}{2}$. Hence we know that $x$ and $a$ lie in the same component of $G\setminus \{b\}$, which is denoted by $G_{abx}$.

Still, we need to show that $G_{abx}$ is open. This directly follows from Lemma \ref{tue-7}
and \cite[Theorem 25.4]{JM}. Hence $G_{abx}$ is a domain in $G\backslash\{b\}$, and so the first statement in the claim is true.\medskip

It follows from the first statement in the lemma that $G_{abx}$ is a subdomain of $G\setminus \{b\}$.
Then the assumptions in the lemma guarantee that the restriction
$f|_{G_{abx}}$ is $\theta$-relative.
For the moment, we need an auxiliary result.

\bas\label{tue-8}
$b\in \partial G_{abx}\subset \partial (G\backslash \{b\})=\partial G\cup\{b\}.$
\eas
Obviously, the relation $b\in \partial G_{abx}$ follows from the rectifiable connectedness of $G$ (cf. \cite[Observation 2.6]{HL}).

For convenience, let $G\backslash \{b\}=G_b$. In the following, we prove the inclusion
$\partial G_{abx}\subset \partial (G_b)$ by contradiction. Suppose on the contrary that there is a point $w\in (\partial G_{abx})\cap G_b$. It follows that $r=\delta_{G_b}(w)>0$, and thus there is a point $u\in \mathbb{B}(w, \frac{1}{1+c}r)$, but $u\notin G_{abx}$.

Since $w\in \partial G_{abx}$, we see that there are $w_n\in G_{abx}$ such that $w_n\to w$ as $n\to \infty$. Without loss of generality, we may assume that all $w_n\in \mathbb{B}(w, \frac{1}{1+c}r)$. Since $u,$ $w_n\in \mathbb{B}(w, \frac{1}{1+c}r)$, we see from Lemma \ref{lem-3.2} \eqref{hwz-131} that there is a curve $\tau'_n\subset \mathbb{B}^{G_b}(w, r)\subset G_b$ connecting $w_n$ and $u$.

On the other hand, the first statement in the claim and \cite[Observation 2.6]{HL} guarantee that $G_{abx}$ is rectifiably connected. Since $a$ and $w_n\in G_{abx}$ (resp. $x$ and $w_n\in G_{abx}$), we know that there is a curve $\tau''_n$ (resp. $\tau'''_n$) in $G_{abx}$ connecting $a$ and $w_n$ (resp. $x$ and $w_n$).
Then $\tau''_n\cup \tau'_n\subset G_b$ (resp. $\tau'''_n\cup \tau'_n\subset G_b$) connects $a$ and $u$ (resp. $x$ and $u$). Hence $u\in G_{abx}$. This contradiction finishes the proof of the assertion.
\medskip

We continue the proof of the claim with the aid of Assertion \ref{tue-8}.
Since
$$\delta_{G}(x)\geq \delta_{G}(z)-|x-z|\geq \delta_{G}(z)-q\delta_{G}(z)> 2q_1\delta_{G}(z)> |x-b|,$$
we infer from Assertion \ref{tue-8} that $$|x-a|=t|x-b|=t\delta_{G_{abx}}(x)<\delta_{G_{abx}}(x),$$ and, again, this assertion guarantees that
\be\label{japan-27} |f(a)-f(x)|\leq \theta(t)\delta_{f(G_{abx})}(f(x))\leq \theta(t)|f(x)-f(b)|.\ee
Hence the proof of the claim is complete.
\medskip

Obviously, $\theta(t)\leq M$. \medskip

Second, we consider the case where $\frac{2c}{1+2c}< t\leq 1.$ In this case, we have
\medskip

\bcl\label{tue-3}
If $\frac{2c}{1+2c}< t\leq 1$ and $0<q_2\leq \frac{1}{(2+c)^2}$, then
$$|f(a)-f(x)|\leq H|f(x)-f(b)|,$$ where $H=H(\theta, c)\geq M$.
\ecl

In this case, by constructing a curve connecting $x$ and $a$ and choosing suitable points from this curve, we shall apply the similar reasoning as in Claim \ref{tue-1} to any two successive points to obtain our desired inequality.

By Lemma \ref{lem-3.2} (\ref{hwz-131}), there is a curve $\beta\subset \mathbb{B}^G(z,(1+c)q_2\delta_{G}(z))$ joining $x$ and $a$ with $${\ell}(\beta)\leq c|x-a|.$$

Define inductively the successive points $x=x_{0},\ldots,x_{k}=a$ of $\beta$ so that $x_{j}$ is the last point of $\beta$ in $\overline{\mathbb{B}}(x_{j-1},(\frac{2c}{1+2c})^{j}|x-b|)$. Obviously, $k\geq 2$, and for $1\leq j\leq k-1$, $$|x_{j-1}-x_{j}|= \left(\frac{2c}{1+2c}\right)^{j}|x-b|$$ and
 $$|x_{k-1}-x_{k}|\leq \left(\frac{2c}{1+2c}\right)^{k}|x-b|.$$

Now, we start to get an upper bound for $k$.
Since for $1\leq j\leq k-1$,
$$l(\beta[x_{j-1},x_{j}])\geq |x_{j-1}-x_{j}|=\Big(\frac{2c}{1+2c}\Big)^{j}|x-b|,$$ we have
$$\sum_{j=1}^{k-1}\left(\frac{2c}{1+2c}\right)^{j}|x-b|\leq l(\beta)\leq c|x-a|= ct|x-b|<c|x-b|.$$
Hence we obtain
\be\label{japan-41} k< \frac{\log 2}{\log(1+2c)-\log (2c)}+1=k_0.\ee

For points $x_1$, $x_0=x$ and $b$, since $|x_1-x_0|=\frac{2c}{1+2c}|x-b|$ and $(1+c)q_2<\frac{1}{2+c}$, the similar method in the proof of \eqref{japan-27} shows that
\be\label{sat-1}|f(x_{1})-f(x_{0})|\leq  M|f(x)-f(b)|,\ee
and for points $x_{0}$, $x_1$ and $x_{2}$, since $|x_2-x_1|\leq \frac{2c}{1+2c}|x_1-x_0|$, similarly, we have
\be\label{sat-2}
|f(x_{2})-f(x_{1})|\leq  M|f(x_1)-f(x_0)|\leq M^2 |f(x)-f(b)|.
\ee

Now, we claim that for $1\leq j\leq k$,
\be\label{sat-3}|f(x_{j})-f(x_{j-1})|\leq  M^{j}|f(x)-f(b)|.\ee

By the relations \eqref{sat-1} and \eqref{sat-2}, we see that to check \eqref{sat-3}, it suffices to discuss the case where $k>2$.
In this case, for $1\leq j\leq k-1$, as in Claim \ref{tue-1}, we use $G_j=G_{x_{j-1}x_jx_{j+1}}$ to denote the component of $G\backslash\{x_j\}$ containing the points $x_{j-1}$ and $x_{j+1}$.
Then by Lemma \ref{tue-7}
and \cite[Theorem 25.4]{JM}, we know that for each $j$, $G_j$ is a domain in $G\backslash\{x_j\}$, and thus the assumptions in the lemma imply that the restriction $f|_{G_{j}}$ is $\theta$-relative.
For points $x_{j+1}$, $x_j$ and $x_{j-1}$, since $$|x_{j}-x_{j+1}|\leq \frac{2c}{1+2c}|x_{j}-x_{j-1}|,$$ it follows from
\beq\label{japan-40}\nonumber
\delta_{G}(x_{j})-|x_{j}-x_{j-1}| &\geq&  \delta_{G}(x_{0})-\sum^{j-1}_{i=0}|x_{i}-x_{i+1}|-|x_{j}-x_{j-1}|
\\ \nonumber &\geq&
\delta_{G}(z)-|x_{0}-z|-\sum^{j-1}_{i=0}|x_{i}-x_{i+1}|-|x_{j}-x_{j-1}|
\\ \nonumber &>& \big(1-(1+4c)q_2\big)\delta_{G}(z)\geq 0\eeq
together with Assertion \ref{tue-8} that for $j\geq 2$,
$$|x_{j}-x_{j+1}|\leq \frac{2c}{1+2c}|x_{j}-x_{j-1}|=\frac{2c}{1+2c}\delta_{D_{j-1}}(x_{j})<\delta_{D_{j-1}}(x_{j}),$$ which shows that
\begin{eqnarray*}
&&|f(x_{3})-f(x_{2})|\leq M\delta_{f(D_1)}(f(x_2))\leq  M|f(x_{2})-f(x_{1})|\leq M^{3}|f(x)-f(b)|, \\
&&\ldots \\
&&|f(x_{k})-f(x_{k-1})|\leq M\delta_{f(D_{k-2})}(f(x_{k-1}))\leq  M^{k}|f(x)-f(b)|.\end{eqnarray*}
Hence the relation \eqref{sat-3} is true.\medskip

It follows from \eqref{sat-3} that
$$
|f(a)-f(x)|\leq  \sum_{j=1}^{k}|f(x_{j})-f(x_{j-1})|
<  kM^{k}|f(x)-f(b)|.$$

Combining the estimate \eqref{japan-41} on $k$, we have
\be\label{japan-28} |f(x)-f(a)|\leq H|f(x)-f(b)|,\ee where
$H=k_{0}M^{k_{0}}\geq M.$ The Claim \ref{tue-3} is proved.\medskip

By Claims \ref{tue-1} and \ref{tue-3}, in fact, we have proved the following conclusion.

\bcor\label{tue-4} The restriction
$f|_{\mathbb{B}^G(z,\frac{1}{(2+c)^2}\delta_G(z))}$ is weakly $H$-quasisymmetric with  $H=k_{0}M^{k_{0}}\geq M$.
\ecor

Finally, we consider the case where $t>1.$

\bcl\label{tue-5}
If $t>1,$ then for any $x$, $a$, $b\in \mathbb{B}^G(z,\frac{1}{(2+c)^3}\delta_G(z))$,
$$|f(a)-f(x)|\leq (1+ct)H^{1+ct}|f(b)-f(x)|.$$\ecl

We shall use Corollary \ref{tue-4} to prove this claim. For any $x$, $a$, $b\in \mathbb{B}^G(z,\frac{1}{(2+c)^3}\delta_G(z))$, Lemma \ref{lem-3.2} (\ref{hwz-131}) implies that
there is a curve $\gamma \subset \mathbb{B}^G(z,\frac{1}{(2+c)^2}\delta_G(z))$ joining $x$ and $a$ with $${\ell}(\gamma)\leq c|x-a|.$$

Define inductively the successive points $x=u_{0}$, $\ldots,$ $u_{\varrho}=a$ of $\gamma$ so that each $u_{i}$ denotes the last point of $\gamma$ in $\overline{\mathbb{B}}(u_{i-1},|b-x|).$ Then
$$\varrho\geq 2, \;\;|u_{i-1}-u_{i}|=|b-x|\;\mbox{ for}\; 1\leq i\leq \varrho-1\;\; \mbox{and}\;\; |u_{\varrho-1}-u_{\varrho}|\leq |b-x|.$$

Now, we come to obtain an upper bound for $\varrho$.
 Since for $1\leq i \leq \varrho-1$, $$\ell(\gamma[u_{i-1},u_{i}])\geq |u_{i-1}-u_{i}|=|b-x|,$$ we have
$$(\varrho-1)|b-x|\leq \ell(\gamma)\leq c|x-a|= ct|b-x|,$$
which implies that
\be\label{fri-2-1} \varrho\leq 1 + ct.\ee

Since $b$ and all $u_i$ $(i\in \{0, \ldots, \varrho\})$ are contained in $\mathbb{B}^G(z,\frac{1}{(2+c)^2}\delta_G(z))$,
Corollary \ref{tue-4} guarantees that
\begin{eqnarray*}
|f(a)-f(x)|&\leq&  \sum_{i=1}^\varrho |f(u_i)-f(u_{i-1})|\leq \varrho H^{\varrho}|f(x)-f(b)|.\end{eqnarray*}
Thus, \eqref{fri-2-1} implies  \be\label{japan-29} |f(a)-f(x)|\leq (1 + ct)H^{(1 + ct)}|f(b)-f(x)|.\ee
The proof of Claim \ref{tue-5} is complete.
\medskip

Now, we are able to finish the proof.
Let  $$\eta'(t)=\begin{cases}
\displaystyle \,\;\;\;\;\;\frac{H}{M}\theta(t),\;\;\;\;\;\;\;
\mbox{if}\;\;0<t\leq \frac{2c}{1+2c},\\
\displaystyle \,\;\;\;\;\; k_0t+k_1,\;\;\;\;\; \,\mbox{if}\;\; \frac{2c}{1+2c}\leq t\leq 1,\\
\displaystyle (1 + ct)H^{(1 + ct)},
\;\mbox{if}\;\; t>1,
\end{cases}
$$ where $$k_0=(1+2c)[(1+c)H^{1+c}-H]\;\;\mbox{and}\;\;k_1=-2c(1+c)H^{1+c}+(1+2c)H.$$
Then it follows from \eqref{japan-27}, \eqref{japan-28} and \eqref{japan-29} that for any $x$, $a$, $b\in \mathbb{B}^G(z,\frac{1}{(2+c)^3}\delta_G(z))$
with $|a-x|=t|b-x|$,
$$|f(a)-f(x)|\leq \eta'(t)|f(b)-f(x)|.$$ This shows that Lemma \ref{lem-4.1} holds.\epf

\subsection{The proof of the implication from \eqref{(1)} to \eqref{(3)} in Theorem \ref{thm-1}}

By using the length metric, we shall show Lemma \ref{lem-4.2} below, from which the proof of
the implication from \eqref{(1)} to \eqref{(3)} in Theorem \ref{thm-1} easily follows. We start with the definition of ``strong geodesic condition" (cf. \cite{BKL}).

\bdefe\label{sat-4}
Suppose $X$ is a rectifiably connected metric space and $G\subsetneq X$ is a domain.
We say that $G$ satisfies the {\it strong geodesic condition} if
every ball $\mathbb{B}(x_0,r)$ of $G$ satisfies the following condition: for every $x\in \mathbb{B}(x_0,r)$ there exists a curve $\gamma$ in $\mathbb{B}(x_0,r)$
connecting $x$ and $x_0$ such that $\ell(\gamma)< r.$
\edefe

The following result is related to the property of the strong geodesic condition for curves.

\bcl\label{cla-4.1} Suppose that $X$ is a $c$-quasiconvex metric space and that both $D\subsetneq X$ and $D'\subsetneq Y$ are domains. If
\ben
\item\label{sat-5}
for any $x_0\in D$ and $0<r< \delta_D(x_0)$, $\mathbb{B}(x_0,r)\subset D$;
\item\label{sat-6}
$D$ satisfies the strong geodesic condition; and
\item
$f:$ $D\to D'$ is $\varphi$-semisolid,
\een
then $f$ is $\theta$-relative with $\theta=\theta_{\varphi, c}$.
\ecl

For the proof, we let $x,$ $y\in D$ with $|x-y|=t\delta_{D}(x)$, where $0< t<1$. We shall find a homeomorphism $\theta:$ $[0,1)\to [0,\infty)$ such that
$$\frac{|f(x)-f(y)|}{\delta_{D'}(f(x))}\leq \theta\Big(\frac{|x-y|}{\delta_{D}(x)}\Big)=\theta(t),$$ where $\theta=\theta_{\varphi, c}$.

To reach this goal, first, we estimate the quasihyperbolic distance between $x$ and $y$, which is as follows:
\be\label{wed-1} k_D(x,y)\leq \frac{1+t}{1-t}.\ee

 Since $y\in \mathbb{B}(x,\frac{1+t}{2}\delta_D(x))\subset D$, we see from the assumptions in the claim that there is a curve $\gamma$ connecting $x$ and $y$ in $\mathbb{B}(x,\frac{1+t}{2}\delta_D(x))$ such that $$\ell(\gamma)< \frac{1+t}{2}\delta_D(x).$$ Hence
$$\delta_D(z)\geq\delta_D(x)-|x-z|\geq \frac{1-t}{2}\delta_D(x)$$  for all $z\in\gamma$,
and as a consequence, we have
$$ k_D(x,y)\leq \int_{\gamma}\frac{|dz|}{\delta_D(z)}\leq \frac{1+t}{1-t},$$
as required.\medskip

Second, by \eqref{vvm-2} in Lemma \ref{lem-3.4}, we have that if $0\leq t\leq \frac{1}{3c}$, then \be\label{wed-2} k_D(x,y)\leq 3ct.\ee

Now, we are ready to construct the needed homeomorphism. Let $\psi(t)=e^t-1$ and
$$\theta_0(t)=\begin{cases}
\displaystyle \;\;\frac{3c(3c+1)}{3c-1}t,\;\;\;\;\;
\mbox{if}\;\;0\leq t\leq \frac{1}{3c},\\
\displaystyle \;\;\;\;\frac{1+t}{1-t},\;\;\;\;\;\;\;\;\;\;\;\; \mbox{if}\;\; \frac{1}{3c}< t< 1.
\end{cases}
$$
Then by \eqref{sun-1} in Lemma \ref{lem-3.4}, we know that
$$\frac{|f(x)-f(y)|}{\delta_{D'}(f(x))}\leq \psi\Big(k_{D'}(f(x),f(y))\Big)\leq \psi\circ\varphi\circ\theta_0(t),$$
and obviously, $\theta=\psi\circ\varphi\circ\theta_0$ is what we want.
Hence Claim \ref{cla-4.1} is complete.\medskip

Next, we are going to show that, in the setting of length metric, the assumptions \eqref{sat-5} and \eqref{sat-6} in Claim \ref{cla-4.1} are naturally satisfied.

\bcl\label{cla-4.3} Suppose that $X$ is a $c$-quasiconvex metric space and that $D\subsetneq X$ is a domain. Then
\ben
\item[\eqref{(1)}]
for any $x_0\in D$ and $0<r< \delta'_D(x_0)$, $\mathbb{B}_d(x_0,r)=\{ x\in X:\; d(x, x_0)<r\}\subset D$; and
\item[\eqref{(2)}]
$D$ satisfies the strong geodesic condition with respect to $d$.
\een
\ecl

Here and in what follows, we always use the notation $\id_{|\cdot|}$ to denote the identity map which is from $(Q,|\cdot|)$ to $(Q,d)$ for any domain $Q$ in $X$.

 Since $(X,d)$ is $\lambda$-quasiconvex for any $\lambda >1$, Lemma \ref{lem-3.6} (\ref{wed-3}) implies that the map $$\id_{|\cdot|}:\;(D,|\cdot|)\to (D,d)$$ is $\lambda$-bilipschtiz, which shows that $(D,d)$ is also a domain.

For any $x_0\in D$, clearly, the ball $\mathbb{B}_d(x_0,r)$ is rectifiably connected, since every two points in $\mathbb{B}_d(x_0,r)$ can be joined by a curve through $x_0$ with length at most $2r$. Then Lemma \ref{lem-3.1} makes sure that for all $r\in (0,\delta'_D(x_0))$, the ball $\mathbb{B}_d(x_0,r)$ is contained in $D$. Hence \eqref{(1)} holds.

To prove \eqref{(2)} in the claim, we let $y$ be an arbitrary point in $\mathbb{B}_d(x_0,r)$. Since $d(x_0,y)<r$, we see that there is a curve $\gamma$ in $X$ joining $y$ and $x_0$ such that
$$\ell_d(\gamma)\leq d(y,x_0)+\frac{1}{2}(r-d(y,x_0))<r.$$
which implies $\gamma\subset \mathbb{B}_d(x_0,r)$, and so, we see that $D$ satisfies the strong geodesic condition with respect to $d$, and hence the proof of Claim \ref{cla-4.3} is complete.\medskip



Suppose that $X$ is a $c$-quasiconvex metric space and that both $G\subsetneq X$ and $G'\subsetneq Y$ are domains. For a homeomorphism $f:$ $G\to G'$, let
$$g=f\circ \id_{|\cdot|}^{-1}:\; (G, d) \to (G',|\cdot|).$$
Then we have the following claim.

\bcl\label{4.4.1}
If $f$ is $\varphi$-semisolid in every subdomain $D$ of $G$, then $g$ is $\theta'$-relative in $D$, where $\theta'=\theta'_{\varphi,c}$.
\ecl

We are going to apply Claims \ref{cla-4.1} and \ref{cla-4.3} to show this claim. We assume that $D$ is an arbitrary subdomain of $G$.
Since for all $x,$ $y\in (D, d)$,
$$k_{f(D)}(g(x),g(y))\leq \varphi(k_D(\id_{|\cdot|}^{-1}(x),\id_{|\cdot|}^{-1}(y)))\leq \varphi(ck'_D(x,y))=\varphi'(k'_D(x,y)),$$ where $\varphi'(t)=\varphi(ct)$ This implies that $g$ is $\varphi'$-semisolid in $D$.
Then it follows from Claim \ref{cla-4.3} that $g$ satisfies all assumptions in Claim \ref{cla-4.1} if the metric $|\cdot|$ is replaced by the one $d$. Claim \ref{cla-4.1} tells us that $g$ is $\theta'$-relative in $D$, where $\theta'=\theta'_{\varphi, c}$. Hence the claim is proved.\medskip

\blem\label{lem-4.2} Suppose that $X$ is a $c$-quasiconvex metric space and that both $G\subsetneq X$ and $G'\subsetneq Y$ are domains. If $f:$ $G\to G'$ is $\varphi$-semisolid in every subdomain of $G$, then \ben
\item\label{wed-8}
there is a homeomorphism $\eta=\eta_{\varphi,c}:$ $[0, \infty)\to [0, \infty)$ such that $f$ is $q$-locally $\eta$-quasisymmetric in $G$, where $q=\frac{1}{(2+c_0)^3c}$;
\item\label{wed-9}
there is a constant $H=H(\varphi,c)$ such that $f$ is $q$-locally weakly $H$-quasisymmetric in $G$.
\een
\elem
\bpf We first prove \eqref{wed-8} in the lemma. Obviously, $(X,d)$ is $c_0$-quasiconvex with $c_0=\frac{\sqrt{3}+1}{2}$.
By applying Claim \ref{4.4.1}, we see that the induced map $$g=f\circ \id_{|\cdot|}^{-1}:\; (G, d) \to (G',|\cdot|)$$ is $\theta$-relative in every subdomain of $G$ with $\theta=\theta_{\varphi, c_0}=\theta_{\varphi}$. So Lemma \ref{lem-4.1} is available now. It follows from this lemma that $g$ is $q'$-locally $\eta$-quasisymmetric with $q'=\frac{1}{(2+c_0)^3}$ and $\eta'=\eta'_{c_0,\theta}=\eta'_{\varphi}$.

Now, we show that $f$ is $q$-locally $\eta$-quasisymmetric, where $q=\frac{1}{(2+c_0)^3c}$ and $\eta(t)=\eta'(ct)$.

For any $x,$ $a,$ $b\in \mathbb{B}^G(z,q\delta_G(z))$ with $|x-a|\leq t|x-b|$, it follows from Lemma \ref{lem-3.6} (\ref{wed-3}) and (\ref{wed-5}) that
$$d(w,z)\leq c|w-z| <cq\delta_G(z)\leq q'\delta_{G}'(z)$$
for all $w\in \mathbb{B}^G(z,q\delta_G(z))$, which implies that $$\mathbb{B}^G(z,q\delta_G(z))\subset \mathbb{B}_d(z,q'\delta_{G}'(z))$$ as a set.
Again using Lemma \ref{lem-3.6} (\ref{wed-3}) we get
$$d(x,a)\leq c|x-a|\leq ct|x-b| \leq ctd(x,b),$$
and then
$$|f(x)-f(a)|=|g\circ \id_{|\cdot|}(x)-g\circ \id_{|\cdot|}(a)|\leq \eta'(ct)|g\circ \id_{|\cdot|}(x)-g\circ \id_{|\cdot|}(b)|=
\eta(t)|f(x)-f(b)|,$$ where $\eta(t)=\eta'(ct)$.
Hence \eqref{wed-8} in the lemma is proved.

The proof of the second statement in this lemma easily follows from \eqref{wed-8} by taking $H=\max\{\eta(1), 1\}$. Hence the proof of Lemma \ref{lem-4.2} is complete.
\epf

\section{The proof of the implication from \eqref{(3)} to \eqref{(1)} in Theorem \ref{thm-1}}\label{sec-5}

This section consists of four subsections. In the first subsection, a mapping property of locally weakly quasisymmetric mappings is given, and the proof of the implication from \eqref{(3)} to \eqref{(1)} in Theorem \ref{thm-1} is contained in the rest three subsections.

\subsection{A mapping property of locally weakly quasisymmetric mappings}

\blem\label{thu-1}
Suppose that $X$ is a $c$-quasiconvex and complete metric space, and that $f:$ $G\to G'$ $q$-locally weakly $H$-quasisymmetric with $H\geq 1$ and $0<q<1$, where $G \subsetneq X$ and $G'\subsetneq Y$ are domains. For $z\in G$, let $Q=\mathbb{B}^{G}(z,r)$, where $0<r\leq\frac{q}{1+2c}\delta_{G}(z)$.
Then $f(\partial Q)=\partial f(Q)$ and $\partial f(\overline{Q})\subset f(\overline{Q})$.
\elem
\bpf
Let $\{f(w_{j})\}_{j=1}^{\infty}$ be a Cauchy sequence in $f({\overline{Q}})$ and $\{w_{j}\}_{j=1}^{\infty}$ the corresponding sequence in $\overline{Q}$.  Then we have
\bcl\label{tue-9}
The sequence $\{w_{j}\}_{j=1}^{\infty}$ is Cauchy.
\ecl
We prove this claim by contradiction. Suppose $\{w_{j}\}_{j=1}^{\infty}$ is not a Cauchy sequence. Then there is $\varepsilon>0$ such that for each positive integer $k$, there is $j(k)>k$ satisfying
$$|w_{k}-w_{j(k)}|\geq \varepsilon.$$ Obviously, there is $z_{k}\in \{w_{k},w_{j(k)}\}$ such that \be\label{japan-39}|z_{k}-w_{1}|\geq \varepsilon/2.\ee
Let $t=\frac{2r}{\varepsilon}$. Then we have $$|w_{1}-w_{j(1)}|\leq 2r= t\varepsilon \leq 2t|z_{k}-w_{1}|.$$
For their images under $f$, we have the following estimate.
\medskip

\bas\label{thu-2} $|f(w_{1})-f(w_{j(1)})|\leq (1+2ct)H^{1+2ct}|f(z_{k})-f(w_{1})|.$ \eas

We divide the proof into two cases. For the first case where $2t\leq 1$, since $|w_{1}-w_{j(1)}|\leq  2t|z_{k}-w_{1}|\leq |z_{k}-w_{1}|$ and
since $w_{1}$, $w_{j(1)}$ and $z_{k}\in \overline{\mathbb{B}}^{G}(z,\frac{q}{1+2c}\delta_{G}(z))\subset \mathbb{B}^{G}(z, q\delta_{G}(z))$, we have
$$|f(w_{1})-f(w_{j(1)})|\leq H|f(z_{k})-f(w_{1})|,$$ as required.

 For the remaining case, that is, $2t> 1$, by Lemma \ref{lem-3.2} \eqref{hwz-131}, there is a curve $\gamma$ in $\mathbb{B}^{G}(z,\frac{1+c}{1+2c}q\delta_{G}(z))$ joining $w_{1}$ and $w_{j(1)}$ such that $${\ell}(\gamma)\leq c|w_{1}-w_{j(1)}|.$$

Define inductively the successive points $w_{1}=a_{0}$, $\ldots,$ $a_{s}=w_{j(1)}$ of $\gamma$ so that for each $i\in \{1, \ldots, s\}$, $a_{i}$ denotes the last point of $\gamma$ in $\overline{\mathbb{B}}(a_{i-1},|z_{k}-w_{1}|).$ Obviously, $s\geq 2$. The following upper bound of $s$ is needed in the proof later on. Since for $i\in \{1, \ldots, s-1\}$,
$$|a_i-a_{i-1}|=|z_k-w_1|\leq \ell(\gamma_{[a_{i-1}, a_i]}),$$
we see that
$$(s-1)|z_k-w_1|\leq \ell(\gamma)\leq c|w_1 -w_{j(1)}|\leq 2ct|z_k-w_1|.$$
Hence \be\label{fri-2}   s\leq 1 + 2ct.\ee

Since all $z_k$ and $a_i$ $(i\in \{0, \ldots, s\})$ are contained in $\mathbb{B}^{G}(z,\frac{1+c}{1+2c}q\delta_{G}(z))\subset \mathbb{B}^{G}(z,q\delta_{G}(z))$,
we know that
\beq\nonumber
|f(w_{1})-f(w_{j(1)})|& \leq & \sum^{s}_{i=1}|f(a_{i})-f(a_{i-1})|\leq  sH^s|f(z_{k})-f(w_{1})| \\ \nonumber
& \leq & (1 + 2ct)H^{1+2ct}|f(z_{k})-f(w_{1})|.
\eeq
Hence the assertion holds.\medskip

Let us proceed with the proof of Claim \ref{thu-2}.
Since $$|z_{k}-w_{1}|\leq 2r= t\varepsilon \leq t|w_{k}-w_{j(k)}|,$$ by replacing $2t$ with $t$, the similar reasoning as in the proof of Assertion \ref{thu-2}
shows that
\be\label{fri-3}\nonumber |f(z_{k})-f(w_{1})|\leq (1+ct)H^{1+ct}|f(w_{k})-f(w_{j(k)})|.\ee

Since $\{f(w_{j})\}_{j=1}^{\infty}$ is a Cauchy sequence, we see that $|f(w_{k})-f(w_{j(k)})|\to 0$ as $k\to \infty$, and so $f(z_{k})\to f(w_{1})$. Then it follows from
Assertion \ref{thu-2} that $f(w_{1})=f(w_{j(1)})$, which is a contradiction since $f$ is homeomorphic. The claim is proved.\medskip

Since $X$ is complete, Claim \ref{tue-9} implies that $\{w_{j}\}_{j=1}^{\infty}$ converges to a point $u\in \overline{Q}$, and so $f(w_{j})\to f(u)\in f(\overline {Q})$,
which shows that the lemma is true.
\epf

\subsection{Locally weak quasisymmetry and ring property}
We start this subsection with the definition of ring property.

\bdefe\label{hwz-7} Suppose $f:$ $G\to G'$ is a homeomorphism, where $G\subsetneq X$ and $G'\subsetneq Y$ are domains. $f$ is said to have the $(M ; \alpha; \beta)$-ring property in $G$ if there are constants $1< \alpha \leq \beta$ and $M>0$ such that for any $B=\mathbb{B}^G(z, r)$ with $z\in G$ and $\beta r<\delta_G(z)$,
$$\diam(f(\overline{B}))\leq M\,\dist\big(f(\overline{B}),\,\partial f(\alpha B)\big),$$
where ${\rm diam}(U)$ (resp. ${\rm dist(U, V)}$) means the diameter (resp. the distance) of the set $U$
(resp. between the two disjoint sets $U$ and $V$).
\edefe


\blem\label{lem-5.1}
Under the assumptions in Theorem \ref{thm-1}, and further, if $f$ is $q$-locally weakly $H$-quasisymmetric in $G$ with $H\geq 1$ and $0<q<1$,
 then $f$ has the $(M ; \alpha; \beta)$-ring property in every subdomain of $G$, where all the constants $M$, $\alpha$ and $\beta$ depend only on $c$, $q$ and $H$.
\elem
\bpf
The assumptions in the lemma imply that
\be\label{japan-38} |f(x)-f(a)| \leq H|f(x)-f(b)|\ee
for any triplet $x$, $a$ and $b$ of points in $\mathbb{B}^{G}(z,q\delta_{G}(z))$ with $|x-a|\leq |x-b|$ and for any $z\in G$.
We shall show that for any subdomain $D$ of $G$, $f$ has the $(M; \alpha; \beta)$-ring property in $D$ with $$M=2H^{2}(H+1),\;\;\alpha =3\;\;\mbox{and}\;\;\beta =6cq^{-1}.$$

For $w\in D$, let $K=\mathbb{B}^D(w,r)$, where $0<r\leq \frac{q}{6(1+c)}\delta_D(w)$.
Then $$\overline{K}\subset\overline{3K}\subset \mathbb{B}^{D}(w,\frac{q}{2c}\delta_{D}(w)).$$
Since Lemma \ref{thu-3} shows $\mathbb{B}^{D}(w,\frac{q}{2c}\delta_{D}(w))\subset\mathbb{B}^{G}(w,q\delta_{G}(w))$, we know that
 $$\overline{K}\subset\overline{3K}\subset \mathbb{B}^{G}(w,q\delta_{G}(w)).$$
We first show a claim.
\medskip

\bcl\label{cla-5.1} $\diam(f(\overline{K}))\leq 2H^2(H+1)\,\dist\big(f(\overline{K}), f(\partial (3K))\big)$.\ecl

Obviously,
it suffices to show that for all $a$, $b$, $y\in \overline{K}$ and $z\in \partial (3K)$,
$$|f(a)-f(b)|\leq 2H^{2}(H+1)|f(y)-f(z)|.$$

Since for any $u\in \partial \overline{K}$, it follows from Lemma \ref{lem-3.3} that
$$\max\{|a-w|, |b-w|\}\leq |u-w|\leq |z-u|.$$ Then we know from \eqref{japan-38} that
$$|f(a)-f(b)|\leq |f(a)-f(w)| + |f(b)-f(w)|\leq 2H|f(u)-f(w)|\leq 2H^2|f(z)-f(u)|.$$
Meanwhile, by Lemma \ref{lem-3.3}, we have
$|y-u|\leq |z-y|$, which implies $$|f(y)-f(u)|\leq H|f(y)-f(z)|,$$
whence
$$|f(z)-f(u)|\leq |f(z)-f(y)|+|f(y)-f(u)|\leq (H+1)|f(y)-f(z)|.$$
Hence
$$|f(a)-f(b)|\leq 2H^2(H+1)|f(y)-f(z)|,$$ as required, and so the claim holds.
\medskip

We are ready to finish the proof of Lemma \ref{lem-5.1}. By Lemma \ref{thu-1} and Claim \ref{cla-5.1}, we see that
$$\diam(f(\overline{K}))\leq 2H^2(H+1)\,\dist\big(f(\overline{K}), f(\partial (3K))\big)
= 2H^2(H+1)\,\dist\big(f(\overline{K}), \partial f( 3K)\big).$$
 Hence the proof of the lemma is complete.
\epf

\subsection{Ring property and relativity}

\blem\label{lem-5.2}
Under the assumptions in Theorem \ref{thm-1}, and further, if $f$ is $q$-locally weakly $H$-quasisymmetric with $H\geq 1$ and $0<q<1$,
then $f$ is $(\theta ; t_{0})$-relative in every subdomain of $G$, where $\theta=\theta_{c,c', M, \alpha, \beta}$ and $t_0=t_0(c, c', M, \alpha, \beta)$.
\elem
\bpf
By Lemma \ref{lem-5.1}, we may assume that $f$ has the $(M ; \alpha; \beta)$-ring property in every subdomain of $G$, where $M=2H^{2}(H+1),$ $\alpha =3$ and $\beta =6cq^{-1}.$
Let $D$ be a subdomain of $G$ and $$t_{0}=\frac{1}{2c(2c\alpha)^{3}\beta}.$$

To prove this lemma, it suffices to prove that
$f$ is $(\theta,t_{0})$-relative in $D$ for some homeomorphism $\theta$: $[0, t_0)\to [0, \infty)$, where $\theta=\theta_{c, c', M, \alpha, \beta}$.
To find such a homeomorphism, we let $x$ and $y\in D$ with
$|x-y|=t\delta_{D}(x)$, where $t\in (0,t_{0})$, and let $m$ be the largest integer with $$2c(2c\alpha)^{m}\beta t< 1.$$
Obviously, $m\geq 3$.
Set $$B_{0}=\mathbb{B}^D(x,t\delta_{D}(x))$$ and   for $1\leq j\leq m$, $$s_{j}=(2c\alpha)^{j}.$$ Then $$2cs_j\beta t\delta_D(x)\leq 2cs_m \beta t\delta_D(x)<\delta_D(x).$$
Let $z_{0}'\in \partial f(D)$ be such that $$|z_{0}'-f(x)|\leq 2\delta_{f(D)}(f(x)).$$ Then there is a curve $\gamma\subset Y$ joining $f(x)$ and $z_{0}'$ such that $${\ell}(\gamma)\leq c'|z_{0}'-f(x)|\leq 2c'\delta_{f(D)}(f(x)).$$

\begin{figure}[htbp]
\begin{center}
\input{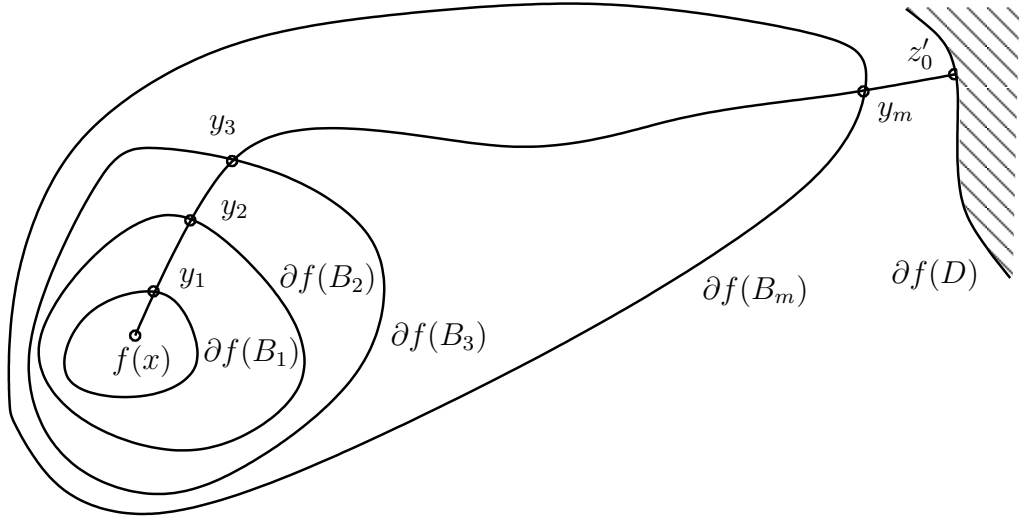}
\caption{The curve $\gamma$ in $Y$ and the related points.}

\label{fqc-fig02}
\end{center}
\end{figure}

Let $$B_j=\mathbb{B}^D(x, s_jt\delta_D(x)).$$ Then it follows that for each $j\leq m,$
\be\label{wed-11}\diam(f(\overline{B}_{j}))\leq M\, \dist( f(\overline{B}_{j}),\partial f(\alpha B_{j}))\leq M\, \dist( f(\overline{B}_{j}),\partial f(2c\alpha B_{j})),\ee and so
\begin{eqnarray*} \dist( f(\overline{B}_{j}), \partial f(D))&\geq& \dist( f(\overline{B}_{m}), \partial f(D))
\geq \dist( f(\overline{B}_{m}),\partial f(2c\alpha B_{m}))\\ \nonumber
&\geq &\frac{1}{M}\diam(f(\overline{B}_{m}))> 0, \end{eqnarray*}
which implies that for each $j\in\{1$, $2$, $\ldots$, $m\}$, the intersection set $\gamma\cap \partial f(B_{j})$ is not empty. Let $y_{j}\in \gamma\cap \partial f(B_{j})$ (see Figure $1$).
Then it follows from \eqref{wed-11} that
 for each $3\leq j\leq m$,
\begin{eqnarray*} |y_{2}-f(x)|&\leq& \diam(f(\overline{B}_{j-1}))
\leq M \dist( f(\overline{B}_{j-1}),\partial f(2c\alpha B_{j-1}))
\\ \nonumber &=& M\, \dist( f(\overline{B}_{j-1}),\partial f(B_{j})).
\end{eqnarray*}

Since for each $j\leq m$, $s_jt<\frac{q}{12c^2}$, we see that $\overline{B}_j\subset \mathbb{B}^D\Big(x, \frac{q}{1+2c}\delta_D(x)\Big)$, and then Lemma \ref{thu-1} implies that
$$\partial f(\overline{B}_{j-1})\subset f(\overline{B}_{j-1}).$$ Hence
$$|y_{2}-f(x)|\leq  M\, \dist(\partial f(\overline{B}_{j-1}),\partial f(B_{j}))\leq
 M\, \dist(\partial f(B_{j-1}),\partial f(B_{j}))
\leq M\, |y_{j}-y_{j-1}|.$$
Summing over the indices $j$, we obtain that
$$
(m-2)|y_{2}-f(x)|\leq  M\,\sum_{j=3}^{m}|y_{j}-y_{j-1}|
 < M{\ell}(\gamma)
\leq 2Mc'\delta_{f(D)}(f(x)).$$

Moreover, since $y\in \overline{\mathbb{B}}(x, t\delta_D(x))$, Lemma \ref{lem-3.2} \eqref{hwz-132} guarantees that $y\in B_1$, and so it follows that
$$|f(x)-f(y)|\leq \diam(f(\overline{B}_{1}))\leq M\dist( f(\overline{B}_{1}),\partial f(B_{2}))\leq M|y_{2}-f(x)|.$$
Hence $$\frac{|f(x)-f(y)|}{\delta_{f(D)}(f(x))}\leq \frac{2M^{2}}{m-2}c'.$$
Since $2c(2c\alpha)^{3}\beta t< 1 \leq 2c(2c\alpha)^{m+1}\beta t$, we have
$$m-2 \geq \frac{\log(1/(2c\beta t))-3\log(2c\alpha)}{\log(2c\alpha)},$$
which implies that $$\frac{|f(x)-f(y)|}{\delta_{f(D)}(f(x))}\leq \frac{2M^{2}c'\log(2c\alpha)}{\log(1/(2c(2c\alpha)^3\beta t))}.$$
Now, let $$\theta(t)= \frac{2M^{2}c'\log(2c\alpha)}{\log(1/(2c(2c\alpha)^3\beta t))}.$$
Obviously, this $\theta$ is the desired. Hence Lemma \ref{lem-5.2} holds.\epf

\subsection{Relativity and semisolidity}

\blem\label{lem-5.3}
Under the assumptions in Theorem \ref{thm-1}, and further, if $f$ is $(\theta ; t_{0})$-relative in every subdomain of $G$, then there exists a homeomorphism $\varphi:$ $[0,\infty) \to [0,\infty)$ such that $f$ is $\varphi$-semisolid in every subdomain of $G$, where $\varphi=\varphi_{c', \theta, t_{0}}$.
\elem
\bpf
Assume that $f$ is $(\theta; t_0)$-relative in every subdomain of $G$.
Then it suffices to prove that $f$ is $\varphi$-semisolid in any subdomain $D$ of $G$, where $\varphi=\varphi_{c', t_0, \theta}$.
To this end, we choose $$t_{1}=\frac{1}{2}\min\Big\{t_{0}, \theta^{-1}\Big(\frac{1}{3c'}\Big)\Big\},$$
 and let $x$ and $y$ be any two points in $D$ with
$k_{D}(x,y)<  t_{1}$. Then it follows from Lemma \ref{lem-3.4} \eqref{vvm-2} that
$$\frac{|x-y|}{\delta_D(x)} \leq 2k_{D}(x,y) < 2t_{1}\leq t_0,$$
which implies
$$|f(x)-f(y)|\leq \theta(2k_{D}(x,y))\delta_{f(D)}(f(x)) \leq \theta(2t_{1})\delta_{f(D)}(f(x)) \leq \frac{1}{3c'}\delta_{f(D)}(f(x)).$$
Again, Lemma \ref{lem-3.4} \eqref{vvm-2} leads to
$$k_{f(D)}(f(x),f(y)) \leq 3c' \frac{|f(x)-f(y)|}{\delta_{f(D)}(f(x))}
\leq  3c'\theta(2k_D(x,y)).$$
Hence, $f$ is $(\psi, t_{1})$-uniformly continuous in the $QH$ metric with $\psi(t)= 3c'\theta(2t)$ (See, for example, \cite{Vai-6} or \cite{Vai-1} for the definition).
Obviously, $\psi(0)=0$.
Since $(D,k_{D})$ is $\lambda$-quasiconvex for any $\lambda>1$, we see from \cite[Lemma 2.5]{Vai-6} or \cite[Lemma 3.2]{Vai-1} that there is a homeomorphism
$\varphi=\varphi_{t_1, \psi}:$ $[0, \infty)\to [0, \infty)$ such that $f$ is $\varphi$-semisolid in $D$. Obviously, $\varphi=\varphi_{c', t_0, \theta}$, and thus Lemma \ref{lem-5.3} is proved.\epf

\section{The proof of the equivalence between \eqref{(1)} and \eqref{(2)} in Theorem \ref{thm-1}}\label{sec-6}

Since the the implication from \eqref{(2)} to \eqref{(1)} is obvious, we only need to show the implication from \eqref{(1)} to \eqref{(2)}.
To prove this implication, it suffices to find constants $K\geq 1$ and $0<\alpha\leq 1$ depending only on $\varphi,$ $c,$ $c'$ such that for any subdomain $D$ of $G$,
\be\label{athm-1} k_{f(D)}(f(x),f(y))\leq K \max\{k_D(x,y)^{\alpha},\;\;k_D(x,y)\}\ee
for all $x,$ $y \in D$. To reach this goal, we need some preparation.

Assume that $D$ is a subdomain of $G$. Then it follows from the first statement in Lemma \ref{lem-4.2} that $f$ is $q$-locally  $\eta$-quasisymmetric with $q=\frac{1}{(2+c_0)^3c}$ and $\eta=\eta_{\theta,c}$, where $c_0=\frac{1+\sqrt{3}}{2}$.

Let
\be \label{thur-1} B_x={\mathbb{B}^D\Big(x,\frac{q}{(2+c)^2}\delta_D(x)\Big)}\;\;\;\; \mbox{and}\;\;\; f_1=f|_{B_x}.\ee
Then \cite[Corollary 3.12]{TV} guarantees that there exist constants $K_0\geq 1$ and $0<\alpha\leq 1$ depending only on $\eta$ such that $f_1$ is $\eta_0$-quasisymmetric, where $$\eta_0(t)=K_0 \max\{t^{\alpha},\;t^{1/\alpha}\}.$$

Next, we show that $f_1$ can be extended to $\overline{B}_x$, and this extension is also $\eta_0$-quasisymmetric. Since $X$ is complete, so is $\overline{B}_x$. Using the facts that the inverse of a quasisymmetric mapping is also quasisymmetric and the image of a Cauchy sequence under a quasisymmetric mapping is also Cauchy, we know that $\overline{f(B_x)}$ is complete. Hence \cite[Theorem 6.12]{Vai-5} makes sure that $f_1$ has a natural extension to $\overline{B}_x$, still denoted by $f_1$, which is also
$\eta_0$-quasisymmetric.

We leave the proof of the inequality (\ref{athm-1}) for a moment and obtain an estimate on the quantity $\frac{|f(x)-f(y)|}{\delta_{f(D)}(f(x))}.$

Let $$t_0=\min\Big\{\varphi^{-1}\Big(\frac{1}{3c'}\Big)\;,\;\frac{q}{(2+c)^4}\Big\}.$$ Obviously, $t_0\leq \frac{1}{81}$. Then for $x$ and $y$ in $D$ with $$|x-y|=t\delta_D(x),$$ we have
\bcl \label{acla-4.3} If $t\leq 2t_0$, then $$\frac{|f(x)-f(y)|}{\delta_{f(D)}(f(x))}\leq c'K_0\Big((2+c)^2q^{-1}\frac{|x-y|}{\delta_D(x)}\Big)^\alpha.$$ \ecl
For the proof, we fix $\varepsilon>0$ and then choose a point $u'\in \partial f(D)$ with
$$|f(x)-u'|\leq \delta_{f(D)}(f(x))+\varepsilon.$$
Since $Y$ is $c'$-quasiconvex, we see that there exists a curve $\gamma$ in $Y$ joining $f(x)$ and $u'$ with
$$\ell(\gamma)\leq c'|f(x)-u'|.$$

Obviously, $\gamma\cap \partial f(B_x)\not=\emptyset$ since $u'\notin f(B_x)$. Furthermore, it follows from Lemma \ref{thu-1} that $\partial f(B_x)=f(\partial B_x)$, and so
$\gamma\cap f(\partial B_x)\not=\emptyset$. Now, we assume that $v'\in \gamma\cap f(\partial B_x)$. Then there is a point $v$ in $\partial B_x$ such that $v'=f(v)$. Hence
$$|f(x)-f(v)|\leq \ell(\gamma)\leq c'|f(x)-u'|\leq c'(\delta_{f(D)}(f(x))+\varepsilon).$$

On the other hand, it follows from Lemma \ref{lem-3.2} \eqref{hwz-132} that $$\overline{\mathbb{B}}\Big(x,\frac{2q}{(2+c)^4}\delta_D(x)\Big)\subset\mathbb{B}\Big(x,\frac{q}{(2+c)^3}\delta_D(x)\Big)\subset B_x,$$ which implies $y\in B_x$. Further, Lemma \ref{lem-3.3} indicates that
$$|x-v|= \frac{q}{(2+c)^2}\delta_D(x).$$ Hence we know from the fact $f_1$ being $\eta_0$-quasisymmetric in $\overline{B}_x$ that
$$\frac{|f(x)-f(y)|}{\delta_{f(D)}(f(x))+\varepsilon}\leq c' \frac{|f(x)-f(y)|}{|f(x)-f(v)|}\leq c'\eta_0\Big(\frac{|x-y|}{|x-v|}\Big)\leq c'\eta_0\Big((2+c)^2q^{-1}\frac{|x-y|}{\delta_D(x)}\Big),$$
from which the claim follows since $(2+c)^2q^{-1}\frac{|x-y|}{\delta_D(x)}< 1$.\medskip

Now, we continue the proof of the inequality (\ref{athm-1}) by dividing the discussions into two cases.

\bca $k_D(x,y)\leq t_0.$\eca

By the choice of $t_0$ and the $\varphi$-semisolidity of $f$ in $D$, we have $$k_{f(D)}(f(x),f(y))\leq \varphi(k_D(x,y))\leq \varphi(t_0)\leq \frac{1}{3c'},$$
and then the inequality \eqref{vvm-2} in Lemma \ref{lem-3.4} gives
\be\label{sun-7} k_{f(D)}(f(x),f(y)) \leq 3c' \frac{|f(x)-f(y)|}{\delta_{f(D)}(f(x))}.
\ee

Since $k_D(x,y)\leq t_0<1$, again, we know from Lemma \ref{lem-3.4} \eqref{vvm-2}  that $$\frac{|x-y|}{\delta_D(x)}\leq 2k_D(x,y)\leq 2t_0.$$
Applying Claim \ref{acla-4.3}, it follows from \eqref{sun-7} that
\be\label{wes-1} k_{f(D)}(f(x),f(y)) \leq 3c'^2K_0\Big((2+c)^2q^{-1}\frac{|x-y|}{\delta_D(x)}\Big)^\alpha<
K_1 k_D(x,y)^{\alpha},\ee
where $K_1=3c'^2K_0(2+c)^{2\alpha}q^{-\alpha}$.

\bca $k_D(x,y)> t_0.$
\eca

By Lemma \ref{lem-3.5}, we see that $(D,k_D)$ is $\lambda$-quasiconvex for all $\lambda>1$. Obviously, the $\varphi$-semisolidity of $f$ in $D$ implies that for any $u,$ $v\in D$, if $k_D(u,v)\leq t_0$, then $k_{f(D)}(f(u),f(v))\leq \varphi(t_0)$. By taking $q=t_0$ and $C=\varphi(t_0)$ in \cite[Lemma 2.3]{Vai-5}, we know that
\be\nonumber k_{f(D)}(f(x),f(y))\leq \lambda\frac{\varphi(t_0)}{t_0}k_D(x,y)+\varphi(t_0)\leq (\lambda+1)\frac{\varphi(t_0)}{t_0}k_D(x,y).\ee
Letting $\lambda\to 1$ yields
\be\label{sun-9} k_{f(D)}(f(x),f(y))\leq K_2 k_D(x,y),\ee where $K_2=\frac{2\varphi(t_0)}{t_0}$.\medskip

Now, we are able to complete the proof of the theorem. By letting $$K=\max\{K_1, K_2\},$$ we see that (\ref{athm-1}) easily follows from \eqref{wes-1} and \eqref{sun-9}, and hence the equivalence between \eqref{(1)} and \eqref{(2)} in Theorem \ref{thm-1} is proved.
\qed

\section{The proof of Theorem \ref{thm-3} and an application}\label{sec-7}

In this section, we prove Theorem \ref{thm-3}.  We begin with the following result showing the invariance of semisolidity under the composition operator of mappings.

\begin{lem}\label{lem-7.1} Suppose $X_i$ are $c_i$-quasiconvex metric spaces, and $G_i\subsetneq X_i$ are domains, where $i=1,$ $2,$ $3$. If $f:$ $G_1 \to G_2$ is  $\varphi_1$-semisolid  in every subdomain of $G_1$ and $g:$ $G_2\to G_3$ is $\varphi_2$-semisolid in every subdomain of $G_2$. Then the composition $g\circ f:$ $G_1\to G_3$ is $\varphi_2\circ\varphi_1$-semisolid in every subdomain of $G_1$.
\end{lem}
\bpf
For any subdomain $D$ in $G_1$, by the hypotheses, we see that $f:$ $D\to f(D)\subsetneq G_2$ is
$\varphi_1$-semisolid and $g:$ $f(D)\to g\circ f(D)\subset G_3$ is $\varphi_2$-semisolid, which shows that for all $x$, $y\in D$,
$$k_{f(D)}(f(x),f(y))\leq \varphi_1(k_D(x,y))$$
and $$k_{g\circ f(D)}(g\circ f(x),g\circ f(y))\leq \varphi_2(k_{f(D)}(f(x),f(y))).$$
 Hence
$$k_{g\circ f(D)}(g\circ f(x),g\circ f(y))\leq \varphi_2 \circ \varphi_1(k_D(x,y)),$$
and so the proof of Lemma \ref{lem-7.1} is complete since the domain $D$ is arbitrarily taken from $G_1$.
\epf

\noindent {\bf The proof of Theorem \ref{thm-3}}:  Under the assumptions, by Theorem \ref{thm-1}, there are homeomorphisms $\varphi_1$ and  $\varphi_2$ from $[0,\infty)$ to $[0,\infty)$ such that
$f$ is $\varphi_1$-semisolid in every subdomain of $G_1$ and $g$ is $\varphi_2$-semisolid in every subdomain of $G_2$, where $\varphi_1$ and $\varphi_2$ depend at most on $c_1$, $ c_2$, $c_3$, $q_1$, $q_2$, $H_1$ and $H_2$. Lemma \ref{lem-7.1} yields that the composition $g\circ f$ is $\varphi_2\circ\varphi_1$-semisolid in every subdomain of $G_1$. Again, it follows from Theorem \ref{thm-1}
that there is a constant $0<q<1$ such that
$g\circ f$ is $q$-locally $\eta$-quasisymmetric in $G_1$, where $\eta:$ $[0,\infty) \to [0,\infty)$ is homeomorphic and depends only on $c_1$, $ c_2$, $c_3$, $q_1$, $q_2$, $H_1$ and $H_2$.
Obviously, $g\circ f$ is $q$-locally weakly $H$-quasisymmetric, where $H=\max\{\eta(1), 1\}$.
Also, it follows from \eqref{hwz-12} in Remark \ref{rem-1.2} that $g\circ f$ is also $H$-quasiconformal.
\qed


\end{document}